\documentclass[11pt, twoside]{article}

\title{Stochastic Localization with Non-Gaussian Tilts and Applications to Tensor Ising Models}
\author{Dan Mikulincer\thanks{Department of Mathematics, University of Washington.}   \and Arianna Piana\thanks{Department of Mathematics, Weizmann Institute of Science.}}

\usepackage{pgf,tikz}
\usetikzlibrary{positioning, arrows.meta}
\usepackage{mathrsfs}
\usepackage[backend=bibtex,bibstyle=numeric, maxnames=10]{biblatex}
\usetikzlibrary{arrows}
\usepackage[a4paper,top=2.5cm,bottom=2.5cm,left=2.5cm,right=2.5cm]{geometry}

\usepackage{polynom}
\usepackage[colorinlistoftodos]{todonotes}
\usepackage{multicol}
\usepackage{mathtools}
\usepackage{amsthm}
\usepackage{amsmath}
\numberwithin{equation}{section}
\usepackage{graphicx}
\usepackage{subcaption}
\usepackage{amssymb}
\usepackage{dsfont}
\usepackage{systeme}
\usepackage{mdframed}
\usepackage{enumitem}
\usepackage{bbm}
\usepackage{bm}
\usepackage{setspace}
\usepackage[colorlinks=true, linkcolor=red, citecolor=cyan]{hyperref}
\theoremstyle{plain}
\usepackage{stmaryrd}
\theoremstyle{plain}
\newtheorem{Thm}{Theorem}[section]

\newtheorem{Lem}[Thm]{Lemma}
\newtheorem{Prop}[Thm]{Proposition}
\newtheorem{Cor}[Thm]{Corollary}

\newtheorem*{theorem*}{Theorem}

\theoremstyle{definition}

\newtheorem*{Que}{Question}

\newmdtheoremenv{theo}{Theorem}
\newtheorem{lemma}[Thm]{Lemma}

\newcommand{\var}{\mathrm{Var}}
\newcommand{\sg}{56}
\newcommand{\spinsg}{200.928}

\newcommand{\op}{\mathrm{op}}
\newcommand{\HS}{_\textup{HS}}
\newcommand{\mix}{_\mathrm{mix}}
\newcommand{\TV}{_\mathrm{TV}}
\newcommand{\dd}{\mathrm{d}}

\newcommand{\Mat}{\mathrm{Mat}}
\newcommand{\Tr}{\mathrm{Tr}}

\newcommand{\EE}{{\mathbb{E}}}
\newcommand{\PP}{{\mathbb{P}}}
\newcommand{\RR}{{\mathbb{R}}}
\newcommand{\SSS}{{\mathbb{S}}}
\newcommand{\ZZ}{{\mathbb{Z}}}
\newcommand{\NN}{{\mathbb{N}}}

\def \DC {\mathcal{C}_n}
\def \vphi {\varphi}
\def \aa {\mathsf{a}}
\def \FF {\mathsf{F}}
\def \aaa {\mathtt{a}}
\def \MMM {\mathtt{M}}
\def \WWW {\mathtt{W}}
\def \CC {\mathsf{C}}
\newcommand{\inj}{\mathrm{inj}}

\newcommand{\GD}{\mathrm{GLD}}

\newcommand{\rank}{\mathrm{Rank}}
\newcommand{\gap}{\mathsf{gap}}

\newcommand{\eps}{\varepsilon}
\newcommand{\poi}{\text{P}}

\newcommand{\dm}[1]{{\color{red}{[[{\bf Dan:} #1]]}}}
\newcommand{\ari}[1]{{\color{magenta}{[[{\bf Ari:} #1]]}}}

\begin{document}
\maketitle

\begin{abstract}
   We present generalizations and modifications of Eldan's Stochastic Localization process, extending it to incorporate non-Gaussian tilts, making it useful for a broader class of measures. As an application, we introduce new processes that enable the decomposition and analysis of non-quadratic potentials on the Boolean hypercube, with a specific focus on quartic polynomials. Using this framework, we derive new spectral gap estimates for tensor Ising models under Glauber dynamics, resulting in rapid mixing.
\end{abstract}

\section{Introduction}
Localization methods are an emerging set of techniques for tackling problems of high-dimensional nature. These ideas originate with the Convex Localization method, simplifying functional and integral quantities by iteratively truncating a given measure along affine hyperplanes \cite{Lovasz1995}. In \cite{RonenThin}, Eldan introduced an adaptation of this approach, known as Stochastic Localization (SL), with several key differences. Rather than truncating the measure, SL perturbs it by introducing a Gaussian factor of the form $\exp( \langle x, \mathsf{d} \rangle - \mathsf{s} \lVert x \lVert_2^2 )$, where the direction $\mathsf{d}$ is chosen randomly. Another distinctive feature is that, while Convex Localization operates in discrete time, SL evolves in continuous time, making its analysis accessible to stochastic analytic tools. Since its original introduction, significant progress has been made in analyzing this method; for more details, see \cite{Chen2021, RonenThin,KlartagLehec, TatLeeVempala}. \par
In recent years, SL has proved to be useful in many other contexts: high-dimensional probability \cite{RonenICM, RonenGaussian, Eldan2018}, mixing for Markov Chains \cite{ChenEldanHit,ChenEldan}, spin glasses \cite{RFIM,Alaoui2022SamplingFT,eldan2022spectral,Montanari_spherical,}, probability on the Boolean hypercube \cite{EG, EldanShamir}, information theory \cite{AhmedMontanari,eldan2020stability}, convex analysis \cite{BoazEli} and sampling \cite{Ghio,MontanariSampling, MontanariSampling2}. The unifying theme across these applications is the following: given an, a priori, complex measure $\mu$ with density $\frac{\dd\mu(x)}{\dd x}= \exp{(-w(x))}$, SL aims to yield better analytic insights into $\mu$ by examining its Gaussian tilts, i.e., $\exp{(-w(x))}\exp{( \langle x, \mathsf{d} \rangle - \mathsf{s} \lVert x \lVert_2^2 )}$. Since the Gaussian distribution is well-understood and possesses numerous desirable analytic and algebraic properties, the hope is that it becomes easier to reason about and analyze the tilted measure. However, a key condition for this approach to succeed is that the potential $w$ must interact with the Gaussian in a way that facilitates analysis. Indeed, most applications therefore focus on measures with some form of convexity, such as when $w$ itself is convex or when dealing with highly structured measures, like quadratic $w$.

In this paper, we extend the stochastic localization framework beyond the Gaussian setting to develop tools applicable to a broader class of measures. Our approach generalizes the localization process to incorporate non-Gaussian tilts, enabling the analysis of potentials that are neither quadratic nor convex. This new framework is versatile and adaptable across various contexts. To showcase its utility, we focus on sampling from tensor Ising models via Glauber dynamics, specifically addressing measures on the Boolean cube where the potential $w$ is a polynomial of degree greater than 2. For clarity, we concentrate on degree 4 potentials as a representative case. In the remainder of the introduction, we will first present the problem of sampling using Glauber dynamics, then explain the limitations of stochastic localization for non-quadratic potentials, and finally outline the necessary adaptations we implement in this work. Our main contributions are summarized below:
\begin{itemize}
    \item We introduce new components and modify existing ones in the stochastic localization process to induce non-Gaussian tilts.
    \item When specializing to degree 4 polynomials, we introduce a variant of the stochastic localization process on tensor spaces and apply our construction to derive new spectral gap estimates for tensor Ising models.
    \item Our bounds also hold, a fortiori, for the spin glass model.
\end{itemize}

\paragraph{Sampling via Glauber Dynamics.} We begin by introducing a motivating question for our generalization of the SL process. Suppose we want to sample from a measure $\mu$ on the $n$-dimensional hypercube $\DC := \{-1,1\}^n$. A common approach is to construct a Markov Chain which converges to $\mu$. 

For measures on $\DC$, Glauber Dynamics (GLD) is a reversible and irreducible Markov chain. At a given time $t \geq 0$, suppose that the current configuration is given by $X_t = x$ where $x = (x_1,\dots,x_n)$.  To define the chain and progress to the next step, we choose a site, or coordinate, $i \in [n]$, and resample $x_i$ according to the conditional distribution of $\mu$ given $(x_j)_{j \neq i}$. In our analysis, we shall consider GLD in continuous time, where times of updates are determined by a homogeneous Poisson process of intensity $n$.

When $\mu$ has full support on $\DC$, it is readily verifiable that the process $X_t$ is ergodic, with $\mu$ as its unique invariant measure. Thus, $\mathrm{Law}(X_t) \xrightarrow{t\to \infty} \mu$, providing a sampling algorithm with easily implementable iterates for $\mu$. To understand the overall efficiency of this algorithm, the key quantity to analyze is the mixing time; in other words, we wish to control $\mathrm{distance}(\mathrm{Law}(X_t),\mu)$ for a given $t \geq 0$.

A well-known method for proving fast mixing for GLD is through its \emph{spectral gap}. Let $L_{\GD}$ denote the generator of GLD, and let $1 = \lambda_1 > \lambda_2 \geq \dots \geq \lambda_{n} \geq -1$ be the eigenvalues of $L_{\GD}$. The spectral gap is defined as $\gap := \lambda_1 - \lambda_2$, and for the semigroup operator $P^t_{\GD}=e^{-tL_{\GD}}$, the mixing time is 
\begin{equation*}
    t\mix(\varepsilon) := \min \big\{ t \geq 0: \ \forall x \in \DC \lVert P_{\GD}^t (x, \cdot) - \mu\rVert\TV \leq \varepsilon \big\}.
\end{equation*} A large spectral gap indicates that the chain mixes quickly, meaning it reaches its stationary in short time. More precisely, these two concepts are related by the following inequality \cite[Theorem 12.4]{LevinPeresWilmer2006}: 
\begin{equation*}
    t\mix(\varepsilon) \leq \left\lceil \frac{n}{\gap} \left( \log\left( \frac{1}{\min_{x \in \DC} \mu(x)}\right) + \log\left( \frac{1}{2\varepsilon} \right) \right) \right\rceil.
\end{equation*}
For GLD, having a spectral gap is equivalent to satisfying a \emph{Poincaré Inequality} (PI). Fix a test function $\varphi: \DC \rightarrow \RR$. The Dirichlet form measures the fluctuations of the function and is defined as
\begin{equation*}
    \mathcal{E}_\mu (\varphi) := \frac{1}{2} \sum_{x, y \in \DC} \big( \varphi(x) - \varphi(y) \big)^2 \mu(x) P_{\GD}(x, y),
\end{equation*}
and the PI is given by
\begin{equation}\label{eq:PIdef}
    \var_{\mu}(\varphi) \leq C_\poi(\mu) \mathcal{E}_\mu (\varphi),
\end{equation}
where the best constant $C_\poi(\mu)$ in \eqref{eq:PIdef}, if it exists, is called the \emph{Poincaré constant} and 
\begin{equation}\label{eq:def_var}
    \var_\mu(\varphi) = \int_{\DC} \varphi^2 \dd \mu - \left( \int_{\DC} \varphi \dd \mu\right)^2 \,.
\end{equation}
The spectral gap is equivalent to the PI: if finite, the Poincaré constant is the reciprocal of the spectral gap. Thus, a larger spectral gap implies a smaller constant in the PI, which means that the variance of functions decays faster over time, leading to faster mixing; for more details, we refer the reader to \cite{LevinPeresWilmer2006, Tetali}.  Below we explain how to obtain bounds on the spectral gap for certain classes of measures using the SL process.
\paragraph{Two questions related to SL.} As explained above, the underlying idea behind SL is to take a measure $\mu$ (on $\RR^n$ or on $\DC$) and tilt it by some Gaussian factor, $\mu(x) \to e^{-Q(x)}\mu(x)$, where $Q$ is a positive semidefinite quadratic form. When the quadratic $Q$ follows some time-evolution in the form of an It\^o process, this transformation makes the measure more tractable for analysis via stochastic calculus. 

To give more details, and explain the connection to mixing times, let us restrict to the discrete case and consider a probability measure $\mu$ on $\DC$. Let $(B_t)_{t \geq 0}$ be a standard Brownian motion on $\RR^n$ with $B_0 = 0$. Given the following (finite) system of stochastic differential equations (SDEs)
\begin{equation}\label{eq:standardSL}
    \dd \FF_t(x)  = \langle x - \aa_t, \CC_t \dd B_t \rangle, \quad \FF_0(x) = 1\,,
\end{equation}
SL is defined as the measure-valued process $\mu_t = \FF_t \mu$, where we now explain the roles of $\aa_t$ and $\CC_t$\footnote{To make the distinctions apparent, for the rest of the paper we use the Serif font to denote the tools for SL.}. First, we see from \eqref{eq:standardSL} that $\FF_t(x)$ is a martingale (since there is no drift term), and therefore, as a probability vector, $\mu = \EE[\mu_t]$, for all $t$, where the expectation is taken over the randomness of the process. 

The vector-valued process $\aa_t$ is chosen to be the barycenter,
\begin{equation*}
     \aa_t = \int_{\DC} x \dd \mu_t(x),
\end{equation*}
and this choice ensures that \(\mu_t\) remains a probability measure. In fact, the total mass is conserved
\begin{equation*}
     \dd \int_{\DC} \mu_t(x) \dd x = \dd \int_{\DC} \FF_t(x) \mu(x) \dd x = \left\langle \int_{\DC} x \FF_t(x) \mu(x) \dd x - \aa_t, \CC_t \dd B_t \right\rangle = 0,
\end{equation*}
and, thanks to the initial condition $\FF_0(x) = 1$, we have that $\mu_t$ is a probability measure at all times. The matrix $\CC_t$ is a matrix-valued process that measures the volatility of the Brownian motion $(B_t)_{t \geq 0}$ and allows to encode various linear constraints on the system.
A simple computation with Itô's formula (see for instance \cite[Fact 14]{ChenEldan}) shows that
\begin{align} \label{eq:SLDef}
    \mu_t(x) &= \exp \left( \int_0 ^t \langle x - \aa_s , \CC_s \dd B_s \rangle - \frac{1}{2} \int_0 ^t \lVert \CC_s (x - \aa_s) \rVert_2 ^2 \dd s \right) \nonumber \\
    &\propto \exp \left( -\frac{1}{2}\langle x, \mathsf{Q}_t x \rangle + \langle \mathsf{L}_t, x \rangle \right)\mu(x), 
\end{align}
where $\mathsf{L}_t$ is some adapted process in $\RR^n$ and $\mathsf{Q}_t = \int_0^t  \CC_s^2 \dd s$ is a matrix-valued process. As $t$ increases, the quadratic term $\mathsf{Q}_t$ increases, which means we are multiplying by Gaussians with smaller and smaller variance. Put differently, the effective support of the measures shrinks, or localizes, over time.
In terms of applicability, this construction raises the following question:
\begin{Que}
    Can we localize a measure by tilting it with a density that is not Gaussian?
\end{Que}
As will become apparent below, there are some rather straightforward ways to induce non-Gaussian localization procedures. However, it will also become apparent that such generalizations are not necessarily useful for potential applications. 
Therefore, it is crucial to develop the new process in a manner that preserves the desirable analytical properties seen in the Gaussian localization while also ensuring its applicability across various non-Gaussian contexts.

In light of this, we will focus on developing the process to establish spectral gaps, as defined in \eqref{eq:PIdef}, a context in which Stochastic Localization has already demonstrated considerable success when the potential is quadratic. Let us now explain the potential usefulness of SL for quadratic functions, which will also make clear that some adaptations are required for more general settings. 

To prove a Poincar\'e inequality via SL, there are several possible approaches, all of which rely on showing a PI for the localized measure $\mu_\tau$ at a suitable finite stopping time $\tau$, i.e. one has to prove
\begin{equation}\label{eq:PItau}
    \var_{\mu_\tau}(\varphi) \leq C_\poi (\mu_\tau) \mathcal{E}_{\mu_\tau}(\varphi)\,,
\end{equation}
where $\varphi:\DC\rightarrow \RR$ is a fixed test function. Most relevant to the present work is the approach adopted in \cite{eldan2022spectral}, whose high-level ideas we now outline. Define $\mathsf{M}_t := \int_{\DC} \varphi \, \dd \mu_t$, to be the running mean of $\vphi$ along the SL process. By the law of total variance, we get
\begin{equation} \label{eq:decompositionVariance}
    \var_\mu(\varphi) = \EE \big[ [\mathsf{M}]_t \big] + \EE \big[\var_{\mu_t}(\varphi) \big],
\end{equation}
where $[\mathsf{M}]_t$ stands for the quadratic variation. 
For matter of intuition, suppose $\mathsf{C}_t$ is chosen in such a way so that $[\mathsf{M}]_t=0$ \emph{almost surely}, then we get a PI:
\begin{equation}\label{eq:PI_EKZ}
   \var_\mu(\varphi) \overset{\eqref{eq:decompositionVariance}}{=}  \EE[\var_{\mu_\tau}(\varphi)] \overset{\eqref{eq:PItau}}{\leq}   C_\poi (\mu_\tau)\EE[\mathcal{E}_{\mu_\tau}(\varphi)]  \overset{(\ast)}{\leq}  C_\poi(\mu_\tau) \EE[\mathcal{E}_{\mu}(\varphi)]\,,
\end{equation}
where in $(\ast)$ we use the easily-verifiable fact that the Dirichlet form associated to $\mu_t$ is a supermartingale, \cite[Proposition 9]{eldan2022spectral} (see also Lemma \ref{lem:dirichlet}). \par
Thus, the argument boils down into carefully choosing the matrix $\CC_t$, which induces a tractable localized measure $\mu_\tau$. These ideas can be instantiated in the following case: Suppose $J$ is an arbitrary symmetric quadratic interaction matrix. For $h \in \RR^n$ an external field, we consider the so-called Ising measure on $\DC$,
\begin{equation} \label{eq:ising_def}
   \mu (x) \propto \exp\big( \langle x, Jx\rangle + \langle h, x \rangle \big).
\end{equation}
For the measure $\mu$, in \cite{eldan2022spectral} the authors choose $\CC_t$ to satisfy two constraints. First, by projecting orthogonally to the level sets of $\varphi$, they force $\EE_{\mu_t}[\vphi]$ to be constant in time, which makes $\var_{\mu_t}(\varphi)$ into a martingale, and hence $\var_{\mu}(\varphi) = \EE\left[\var_{\mu_t}(\varphi)\right]$. Second, by making sure that $\mathrm{Image}(\CC_t)$ is always contained in an appropriate subspace, they ensure that $J_t = J - \int_0^t \CC_s^2 \dd s$, which appears in \eqref{eq:SLDef}, is strictly decreasing as long as $\rank(J_t) \geq 2$. If $\tau := \min\{ t : \rank(J_t) = 1 \}$, then $\mu_\tau$ corresponds to a rank-$1$ matrix $J_\tau$. In other words, the process decomposes the measure $\mu$ into a mixture of simplified measures of the form
\begin{equation}\label{eq:EKZ_mixture}
    w_{u,v}(x)\propto \exp \big(\langle u, x\rangle^2+\langle v, x\rangle \big).
\end{equation}
Moreover, the integral of the test function $\varphi$ is preserved, i.e. $\int \varphi \dd w_{u,v} = \int \varphi \dd \mu$ a.s., and $\|u\|^2 = \|J_t\|_{\mathrm{op}} \leq \|J\|_{\mathrm{op}}$ . Such measures are simpler to treat and satisfy a PI, with $C_\poi(w_{u,v}) \leq (1-\lVert u\lVert_2^2)^{-1}$. Together with \eqref{eq:PI_EKZ} this implies $C_\poi(\mu) \leq (1-\|J\|_{\mathrm{op}})^{-1}$.

The fact that $J$ is a quadratic form plays a crucial role in all of the above steps. It is this fact that forces $J_t$ to remain a matrix at all times and to eventually localize to a rank $1$ matrix. Clearly, if $J$ were any other type of function, it would not make sense to ask whether $J - \int_0^t \CC_s^2 \dd s$ is decreasing. Thus to handle more general potentials $J$ we need to alter the localization scheme, which leads us to the following question.
\begin{Que}
    Can similar techniques be applied to interactions of higher degrees?
\end{Que}
From now on we shall focus on answering this question, by considering one of the simplest non-quadratic models of interactions. Specifically, we explain how to tackle degree $4$ polynomials. I.e., we replace the matrix $J$ with a fourth-order tensor. 
\paragraph{Tensor Ising models and a na\"ive approach for fourth-order tensors.}
Let $T\in (\mathbb{R}^n)^{\otimes 4}$ be a fourth-order tensor and for $x \in \RR^n$ write $T(x) = T(x,x,x,x) = \langle T, x^{\otimes 4}\rangle\HS$ when we regard $T$ as a multi-linear operator. Equivalently, $T(x) = \left\langle x^{\otimes 2}, Tx^{\otimes 2} \right\rangle\HS$ if instead we think about $T$ as an linear operator acting on matrices, see Section \ref{sec:tensors} for more details.

With the tensor $T$, consider the measure $\nu(x) \propto \exp(T(x))$ on $\DC$. The measure $\nu$ is also called a \emph{tensor Ising model}. It is similar to the classical Ising model on the vertex set $[n]$, except that we allow higher-order interactions, in this case $4$, between vertices.  When dealing with tensors, a seemingly natural approach to extend stochastic localization is to define the following version of SL: $\nu_t = \mathtt{F}_t \nu$, where $\mathtt{F}_t$ evolves according to the SDEs
\begin{equation} \label{eq:NaiveSL}
    \dd \mathtt{F}_t(x) = \langle x^{\otimes 2} - \aaa_t, \MMM_t \dd \WWW_t \rangle\HS, \quad \mathtt{F}_0(x) = 1,
\end{equation}
with
\begin{equation*}
    \aaa_t = \int_{\DC} x^{\otimes 2} \dd \nu_t(x),
\end{equation*}
and where $\WWW_t$ is a Dyson Brownian motion\footnote{With a slight abuse of terminology, by Dyson Brownian motion, we refer to a symmetric $n^2\times n^2$-matrix where each entry is as an independent standard Brownian motion.}, and $\MMM_t$ is a fourth-order tensor operating on matrices. The objective is to tilt the measure using a fourth-order tensor.
By applying Itô's formula, we obtain the following expression:
\begin{align}
    \nu_t(x) &= \exp \left( \int_0^t \langle x^{\otimes 2} - \aaa_s, \MMM_s \dd \WWW_s \rangle\HS - \frac{1}{2} \int_0^t \lVert \MMM_s(x^{\otimes 2} - \aaa_s) \rVert^2 \HS\dd s \right) \mu(x) \nonumber\\
    &\propto \exp\left( \left\langle x^{\otimes 2},\left(T -\frac{1}{2}\int_0^t \MMM_s^2 \dd s\right)x^{\otimes 2}\right\rangle\HS + \left\langle \mathtt{L}_t, x^{\otimes 2} \right\rangle\HS\right), \label{eq:problems}
\end{align}
where 
\begin{equation} \label{eq:toobig}
    \mathtt{L}_t= \int_0^t \big(\MMM_s\dd \WWW_s+\MMM_s\aaa_s\dd s\big)\,.
\end{equation}
Thus, at a superficial level, by an appropriate choice of $\mathtt{M}_t$ we can tilt the measure by any desirable fourth-order tensor. In particular, by following the steps of \cite{eldan2022spectral} we can apparently simplify the arbitrary tensor $T$ into a lower-rank tensor. 

However, as we now explain, this na\"ive localization procedure is unsuitable for applications.
First, in the standard SL method \eqref{eq:EKZ_mixture}, we encountered a random linear term $\langle h,x\rangle$, which results in tilting by a \emph{product measure}. In the context of GLD for the Ising model, which updates according to the conditional distributions, as in \cite{eldan2022spectral}, tilting by a product measures is manageable, and has no adverse effect on the spectral gap. However, the last term in \eqref{eq:problems} is quadratic in $x$, and since $\mathtt{L}_t$ is random, it can potentially become large. In fact, for certain choices of $\mathtt{L}_t$, it can be shown that the spectral gap deteriorates significantly. More generally, the term $\langle \mathtt{L}_t, x^{\otimes 2}\rangle_{\mathrm{HS}}$ is not convex in general, unlike the linear term in \eqref{eq:SLDef}, which is prohibitive for a large range of applications.

Furthermore, at a technical level, we treat $T$ and $\mathtt{M}_t$ as operators on matrices, i.e. elements of $\mathrm{Mat}(\RR^{n^2},\RR^{n^2})$, forgetting some of their tensor structure. Thus the term $T(x) - \frac{1}{2}\int_0^t \MMM_s^2(x) \dd s$ collapses to a rank-$1$ element in the matrix space $\mathrm{Mat}(\RR^{n^2},\RR^{n^2})$, rather than retaining the full structure of a tensor. 
Concretely, a symmetric rank-$1$ element in $\mathrm{Mat}(\RR^{n^2},\RR^{n^2})$ is of the form $M \otimes M$ for some $n\times n$ matrix $M$, rather than $v^{\otimes 4}$ for a vector $v\in \RR^n$.

\paragraph{Contributions.} To address the challenges mentioned above, we devise a new decomposition scheme that we apply in several steps. We begin by constructing a variant of the SL method which essentially, up to an arbitrarily small error, decomposes the measure $\dd \mu(x) \propto \exp(T(x))$ into a family of \emph{non-negative} measures of the form 
\begin{equation*}
\mu_{M,R}(x) \propto \exp \left( \langle x, M x \rangle ^2+\langle x, Rx\rangle\right),
\end{equation*}
such that 
\begin{equation*}
\mu = \int \mu_{M,R} \eta(\dd M, \dd R).
\end{equation*}
In the expression above, note that we completely omitted the term $\langle \mathtt{L}_t,x^{\otimes 2}\rangle\HS$. In exchange, we can no longer guarantee that $\mu$ decomposes into a mixture of probability measure, as $\mu_{M,R}$ is not necessarily normalized to have total mass $1$. The key point however is that, as we shall show, the norms of $M$ and $R$ are bounded almost surely in terms of $\|T\|_{\mathrm{inj}}$, the \emph{injective norm} of $T$, see Section \ref{sec:tensors} for the definition. This property is still maintained after re-normalizing $\mu_{M,R}$ to be a probability measure, and will eventually allow us to establish a spectral gap. To understand how a spectral gap behaves for general non-negative measures the following quantity will play a critical role: If $\nu$ is any non-negative measure on $\DC$ with total mass $\EE_{\nu}[1]$, and $\vphi:\DC\to\RR$, we define the \emph{normalized variance} of $\vphi$ with respect to $\nu$ by,
$$\overline{\var}(\vphi):=\EE_\nu\left[\vphi^2\right] - \frac{\EE_\nu\left[\vphi\right]^2}{\EE_\nu\left[1\right]}.$$
Above $\EE_\nu\left[\vphi^2\right] = \sum\limits_{x\in \DC} \vphi(x)^2\nu(x)$  and so it is straightforward to verify (see also Lemma \ref{lem:nonnegativesg}) that $\overline{\var}_{\nu}(\vphi) = \EE_\nu\left[1\right]\var_{\tilde{\nu}}(\vphi)$ for the normalized measure $\tilde{\nu} = \frac{\nu}{\EE_\nu\left[1\right]},$ which justifies the name of $\overline{\var}$. As it turns out, when working with general measures it is desirable to preserve the normalized variance, rather than the standard variance, which is how we design our decomposition.

To achieve the decomposition we change the SL process in two distinct ways:
\begin{itemize}
    \item We first remove the barycenter $\aa_t$ from the construction in \eqref{eq:NaiveSL}. For this reason the total mass of $\mu_t$ is going to change along the SL process, which is why we will have to deal with general non-negative measures. The constraint in the matrix $\MMM_t$ is then adapted to deal with such measures and to preserve the normalized variance.
    \item After removing $\aa_t$, we introduce a new process $v_t$ to take its place. We choose the new process to ensure that the term $\mathtt{L}_t$ in \eqref{eq:toobig} is arbitrarily small.

\end{itemize}
The above decomposition constitutes a rank-$1$ decomposition of matrices, rather than a decomposition of tensors. To obtain a full tensor decomposition, we iteratively refine our approach by decomposing each matrix term separately. Specifically, applying our decomposition to one of the copies $M$ yields an expression of the form:
\begin{equation*}
\Tilde{\mu}(x) \propto \exp \left( \langle u, x \rangle^2\langle x, M x \rangle+\langle x, Rx\rangle\right).
\end{equation*}
Here, we decompose one copy of $M$ into the rank-$1$ matrix $u \otimes u$. We further apply our SL variant to the other copy of $\langle x, M x \rangle$. To complete the decomposition, we apply SL as in \cite{eldan2022spectral} to the quadratic term $\langle x, Rx\rangle$.
Figure \ref{fig:schema} illustrates this iterative refinement process when applied to a tensor $T$.
This results in the following decomposition, which also incorporates certain lower-order terms.
\begin{figure}[h!]
    \centering
    \resizebox{7.5cm}{!}{\begin{tikzpicture}[
    arrow/.style={-Stealth, thick},
    line/.style={draw, thick},
    dottedline/.style={draw, thick, dotted}
]

\node (A) at (0.6, 0) {\( \exp\big( T(x)\big) \)};
\node (B0) at (-2.1, -2) {$\exp\big($};
\node (B1) at (-0.9, -2) {$ \langle x, M x \rangle$};
\node (B2) at (0.4, -2) {$  \langle x, M x \rangle $};
\node (B4) at (2.4, -2) {$+\quad \langle x, Rx\rangle\big)$};
\node (B3) at (0, -2) {};
\node (C0) at (-2.1, -4) {$\exp\big($};
\node (C1) at (-0.9, -4) {$  \langle x, M x \rangle $};
\node (C2) at (0.4, -4) {$  \langle v,  x \rangle^2 $};
\node (C4) at (2.4, -4) {$+\quad\langle x, Rx\rangle\big)$};
\node (D0) at (-2.1, -6) {$\exp\big($};
\node (D1) at (-0.9, -6) {$  \langle u,  x \rangle^2 $};
\node (D2) at (0.4, -6) {$  \langle v,  x \rangle^2 $};
\node (D4) at (2.4, -6) {$+\quad\langle x, Rx\rangle\big)$};
\node (E0) at (-2.1, -8) {$\exp\big($};
\node (E1) at (-0.9, -8) {$  \langle u,  x \rangle^2 $};
\node (E2) at (0.4, -8) {$  \langle v,  x \rangle^2 $};
\node (E4) at (3, -8) {$+\quad\langle x, w\rangle^2 +\langle \ell,x\rangle\big)$};

\draw[thick] (B2) --node[right]{Thm \ref{thm:seconddecomp}}++ (C2);
\draw[dottedline] (B1) -- (C1);
\draw[thick] (C1) --node[left]{Thm \ref{thm:seconddecomp}}++ (D1);
\draw[dottedline] (C2) -- (D2);
\draw[dottedline] (D1) -- (E1);
\draw[dottedline] (D2) -- (E2);
\draw[thick] (D4) --node[right]{Thm \ref{thm:rank1decompo}}++ (E4);
\draw[thick] (0.6,-0.2) --node[right]{Thm \ref{thm:firstdecomp}}++ (0,-1.5);
\draw[dottedline] (B4) -- (C4);
\draw[dottedline] (C4) -- (D4);
\end{tikzpicture}}
    \caption{Decomposition of a fourth-order tensor.}
    \label{fig:schema}
\end{figure}

\begin{Thm} \label{thm:maindecompo}
    Let $\varphi: \DC \rightarrow \RR$ be a test function and let $\mu$ be a measure on $\DC$ given by
    \begin{equation*}
		\mu(x) \propto \exp(T(x)),
	\end{equation*}
	where $T$ is a positive definite symmetric fourth-order tensor with zero diagonal entries. Then, for every $\delta>0$, there exists a decomposition of $\mu$ as
	\begin{equation*}
		\mu = \int \mu_{u, v, w, \ell,\psi} \eta(\dd u, \dd v, \dd w, \dd\ell,\dd\psi),
	\end{equation*}
	where with $u, v, w,\ell \in \RR^n$ and $\psi:\DC\to \RR$, and $\mu_{u, v, w, \ell,\psi}$ is a non-negative measure of the form
	\begin{equation*}
		\mu_{u, v, w, \ell,\psi}(x) \propto \exp\left(  \langle u, x \rangle^2 \langle v, x \rangle^2 + \langle x, w \rangle^2  + \langle \ell ,x\rangle +\psi(x)\right).
	\end{equation*}
 Moreover, we have the following properties:
 \begin{itemize}
     \item $\| u\|_2^2, \| v \|_2^2 \leq 2\sqrt{\lVert T \lVert _\inj} $, $\| w\|_2^2 \leq 4n \lVert T \lVert_\inj$,  $\eta$-almost surely. 
     \item $\max\limits_{x\in \DC}|\psi(x)|\leq\delta$,  $\eta$-almost surely. 
     \item The variance is decomposed as follows 
     $\var_\mu(\varphi) \leq \int \overline\var_{\mu_{u, v, w, \ell,\psi}}(\vphi) \eta(\dd u, \dd v, \dd w, \dd\ell,\dd\psi) +\delta$.
 \end{itemize}
 \end{Thm}
As can be seen, extending from quadratic to quartic potentials involves several non-trivial modifications, which is why we chose to focus on this case. However, we note that the modifications we described are applicable in more general settings and can be used to induce a variety of other decompositions. Furthermore, as will be seen from the proof, the discrete cube $\DC$ does not play a central role in the decomposition and could replaced, for example, with the sphere $\mathbb{S}^{n-1}$ or any of its measurable subsets.

With Theorem \ref{thm:maindecompo} we proceed with the same strategy from \cite{eldan2022spectral}, as outlined above, and obtain spectral gaps for tensor Ising models. 
\begin{Thm} \label{thm:tensorgap}
		Let $\mu$ be a measure on $\DC$, given by,
		$$\mu(x) \propto \exp(T(x)),$$
		for some positive definite degree-$4$ tensor $T$. Assume that $\|T\|_{\mathrm{inj}}\leq \frac{1}{\sg n}$. Then,
		$$C_{\textnormal{\poi}}(\mu)\leq \frac{1}{1 - \sg n\|T\|_{\mathrm{inj}} }\,.$$
	\end{Thm}

Let us make several remarks concerning Theorem \ref{thm:tensorgap}.  First, at this level of generality, the requirement $\|T\|_{\mathrm{inj}}=O(\frac{1}{n})$ is necessary for rapid mixing. This condition arises naturally when working in $\DC$, where each element has norm $\sqrt{n}$. Consequently, we recover the correct order of magnitude for the bound on $\|T\|_{\mathrm{inj}}$. However, the constant $\frac{1}{\sg}$ is likely suboptimal. This sub-optimality stems from two main factors. 
The main reason for sub-optimality lies in our use of \emph{Dobrushin's condition}, see \eqref{eq:dobru}, to bound the spectral gap for rank-one tensor Ising models. While the condition is tight for rank-one matrices, it leaves a gap for higher-order tensors. Tightening this condition for tensors seems feasible, but will likely require some new ideas. The other reason for the sub-optimality is technical and seems to be an artifact of our method. In the decompositions discussed above, we cannot ensure that the matrix $M$ is positive semi-definite, a required condition for our method. We then need to shift $M$ to become positive semi-definite, which artificially increases its operator norm and induces the lower order terms. This shift also causes the sub-optimality to further degrade when dealing with higher degrees.

To understand the extent of the sub-optimality in our constant $\frac{1}{\sg}$ we focus on a specific, simple, tensor Ising model sometimes called tensor Curie-Weiss. 
\begin{lemma} \label{lem:4CW}
Let $T(x) = \frac{1}{n^3}(\sum x_i)^4$, and for $\beta>0$ consider the tensor Curie-Weiss measure
$$\mu_\beta(x) \propto \exp\left(\beta T(x)\right).$$
Suppose that $\beta \gtrsim 0.50425$. Then, the mixing time of Glauber dynamics, associated with $\mu_\beta$, is at least exponential. 
\end{lemma}
We prove Lemma \ref{lem:4CW} in Section \ref{sec:CW} as a consequence of a more general result for arbitrary degree tensors. We can see that the tensor $T$ in the Lemma is a degree-$4$ tensor of rank one. Moreover, since $\|\beta T(x)\|_{\mathrm{inj}} = \frac{\beta}{n},$ we can see that as $\beta$ varies from $\frac{1}{\sg}$ to $0.50425$ the measures $\mu_\beta$ undergo a transition from rapid mixing of GLD to slow mixing. We believe the tensor Curie-Weiss is the worst-case example for GLD, which shows that our bound is probably off by about a factor of $28$. Still, to our knowledge, there is no other general bound with an explicit constant in this setting.

A few days after we uploaded our paper to the arXiv, the authors of \cite{samanta2024mixing} uploaded their independent work about mixing times for GLD in the tensor Curie-Weiss model. The paper includes a more systematic treatment of this specific model and in particular a different proof of Lemma \ref{lem:4CW} and Proposition \ref{prop:TensorCW} with improved lower bounds.\\

For another application, we specialize Theorem \ref{thm:tensorgap}, which is a general statement about tensor Ising models, to specific cases of interest. Namely, we consider the degree $4$ spin-glass model, in which the tensor $T$ has independent Gaussian entries.
\begin{Cor} \label{cor:spinglass}
	Let $T$ be a symmetric fourth-order tensor with    independent $\mathcal{N}(0,\frac{1}{n^3})$ off-diagonal entries and $0$ everywhere else. For fixed $\beta > 0$, let $\mu$ be a probability measure on $\DC$ given by,
		$$\mu_\beta(x) \propto \exp\left(\beta T(x)\right).$$
		If $\beta \lesssim\frac{1}{\spinsg}$, then
		$$C_\textnormal{\poi}(\mu) \leq \frac{1}{1- \spinsg\beta}.$$
	\end{Cor}
Corollary \ref{cor:spinglass} follows from Theorem \ref{thm:tensorgap} using an appropriate bound on the injective norm of a Gaussian tensor. Again, while the constant $\frac{1}{\spinsg}$ is not optimal, it remains, to the best of our knowledge, the best-known bound in this setting, compared to the implicit constant in \cite{adhikari2024spectral}.

Finally, we note that our method extends to certain tensor models of higher degrees. However, the rank of the decomposition grows rapidly with the degree, causing a sharp decline in the spectral gap guarantees. These issues are further discussed in Section  \ref{sec:higher}.
\paragraph{Further related works.}
As mentioned above, our decomposition result in Theorem \ref{thm:maindecompo} should be viewed as an extension of the rank-$1$ decomposition theorem in \cite[Theorem 12]{eldan2022spectral}. For completeness, we state the previous theorem and then discuss the important differences.

\begin{Thm} (\cite[Theorem 12]{eldan2022spectral}) \label{thm:rank1decompo}
    Let $\varphi: \DC \rightarrow \RR$ be a test function and let $\mu$ be a measure on $\DC$ given by
    \begin{equation*}
		\mu(x) \propto \exp(\langle x,Jx\rangle),
	\end{equation*}
	where $J$ is a positive definite matrix. Then, there exists a decomposition of $\mu$ as
	\begin{equation*}
		\mu = \int \mu_{u,\ell} \eta(\dd  u, \dd \ell),
	\end{equation*}
	where $u,\ell \in \RR^n$, and $\mu_{u, \ell}$ is a probability measure of the form
	\begin{equation*}
		\mu_{u, \ell}(x)\propto \exp \big(\langle u, x\rangle^2+\langle \ell, x\rangle \big).
	\end{equation*}
 Moreover, we have the following properties:
 \begin{itemize}
     \item $\|u\|_2^2 \leq \|J\|_{\mathrm{op}}$,  $\eta$-almost surely.
     \item The variance is decomposed as follows $\var_\mu(\varphi) = \int \var_{\mu_{u,\ell}}(\varphi) \eta(\dd u, \dd \ell)$.
 \end{itemize}
 \end{Thm}
As can immediately be seen, Theorem \ref{thm:rank1decompo} decomposes $J$ into rank-$1$ matrices $u\otimes u$. In contrast, since we encode more linear constraints in the tensorized stochastic localization process, Theorem \ref{thm:maindecompo} can only offer a rank-$8$ decomposition, augmented with some lower order terms, which leads to the suboptimality of our spectral gap estimates. The upshot is that the `linear' term $\mathtt{L}_t$ from \eqref{eq:problems} does not appear in Theorem \ref{thm:maindecompo}. As discussed above, omitting this term is crucial for applications.

Concerning spectral gaps, other than \cite{eldan2022spectral} the SL process was instrumental in establishing several other spectral gap results. In \cite{AhmedMeanfield, Alaoui2022SamplingFT}, the authors discretize the SL process and produce a new sampling algorithm for the Sherington-Kirkpatrick (SK) and spin-glass models, using AMP iterations. This new algorithm affords non-trivial sampling guarantees for the SK model when the inverse temperature satisfies $\beta < 1/2$. The paper \cite{Celentano} later extended this result to the conjectured optimal threshold $\beta < 1$. The work in \cite{AKV} further improves the fast mixing results for the GLD in the Ising model, implying fast mixing for the SK-model up to $\beta \approx 0.295$. SL has also been applied to other models, including the random field Ising model \cite{RFIM}, the spherical spin-glass \cite{Montanari_spherical, huang2024weak}, and the sparse random Ising model \cite{LMRW}. For a survey of recent developments and applications to machine learning, see \cite{MontanariSampling}.

To compare with our Theorem \ref{thm:tensorgap}, perhaps the most relevant result is from \cite{AJKTV}. In that work, by analyzing the spectral independence condition, the authors show that for \emph{any} measure $e^{-H(x)}$ on  $\DC$, the spectral gap is bounded as long as $\max_{x\in \DC}\|\nabla^2H(x)\|_{\mathrm{op}}$ is sufficiently small. 
For a degree-$4$ symmetric tensor, as in Theorem \ref{thm:tensorgap}, we have
 $$\max_{x\in \DC}\|\nabla^2T(x)\|_{\mathrm{op}} = 12\max_{x\in \DC} \|T(x,x,\cdot,\cdot)\|_{\mathrm{op}} = 12n\|T\|_{\mathrm{inj}}.$$
Thus, \cite{AJKTV} implies the existence implicit of an implicit small constant $B$, such that when $\|T\|_{\mathrm{inj}} < \frac{B}{n}$, GLD exhibits rapid mixing. Our result improves upon this bound by providing an explicit constant $\frac{1}{\sg}$, replacing $B$.
However, for higher-degree tensors, the estimates in \cite{AJKTV} scale more favorably compared to those in Section \ref{sec:higher}. Despite this, our approach offers several new insights, beyond providing an explicit constant. First, as mentioned above, our decomposition extends to continuous models in the sphere, whereas the spectral independence condition is inherently discrete. More importantly, Theorem \ref{thm:maindecompo} establishes a connection between general tensor Ising models and low-rank models, enabling the transfer of analytic insights from the latter to the former. This decomposition has potential applications beyond spectral gap analysis. For instance, if it could be shown that low-rank tensor models satisfy a modified log-Sobolev inequality, as in \cite{anari2022entropic}, the same result would extend to general models through our framework. Establishing such an inequality, which would imply near-linear mixing time for GLD, is an interesting question that is beyond the scope of this work.

\paragraph{Overview.} The structure of the paper is outlined as follows. In Section \ref{sec:tensors}, we review some details about fourth-order tensors and establish the necessary properties. Section \ref{sec:TSL} is dedicated to constructing the tensorized Stochastic Localization process. In Section \ref{sec:spectral_gap}, we apply our process to establish estimates on spectral gaps. Section \ref{sec:CW} is dedicated to the tensor Curie-Weiss model and in Section \ref{sec:higher} we discuss further generalizations for higher-order tensors.

\paragraph{Conventions.} We denote the space of \( m_1 \times m_2 \) matrices over \( \mathbb{R}^n \) by \( \Mat(\mathbb{R}^{m_1}, \mathbb{R}^{m_2}) \). Let \( A, B \in \Mat(\mathbb{R}^m, \mathbb{R}^m) \) for some \( m \in \mathbb{N} \). The Hilbert-Schmidt inner product is defined as $\langle A, B \rangle{\HS} := \Tr(A^* B).$ Moreover, the operator norm is defined as $\lVert A\lVert_\op := \sup \{\lVert Av\lVert_2: \lVert v\lVert_2\leq 1 \text{ and } v \in \RR^m \}$. Let $x,y \in \RR^n$, we define the Euclidean scalar product as $\langle x, y\rangle = x^T y$, the tensor product as $x\otimes y=xy^T$ and the $m$-th tensor power as $x^{\otimes m}= x \otimes \dots \otimes x.\\$
Let $f_1, f_2: \Omega \rightarrow \mathbb{R}$ be two non-negative functions over the same space $\Omega$. We write $f_1(x) \propto f_2(x)$ if there exists a positive constant $Z$ such that $f_1(x)=\frac{1}{Z} f_2(x)$ for every $x \in \Omega$.\\
Given a measure $\mu$ on $\DC$, we say it admits a decomposition if there exists a measure $m$ on some index set $\Theta$ and a family of probability measures $(\mu_\theta)_{\theta \in \Theta}$ on $\DC$ such that $\mu$ can be rewritten as $\mu(A) = \int_{\Theta}\mu_\theta(A) \dd m(\theta)$ for all $A \subset \DC $ measurable. We call $m$ the mixing measure.
With a slight abuse of notation, we shall also identify $\mu$ with its probability function, so that for $x\in \DC$ we write $\mu(x)$ instead of $\mu(\{x\})$.

\paragraph{Acknowledgments.} A.P. is supported by the Israeli Council for Higher Education (CHE) via the Weizmann Data Science Research Center, and by a research grant from the Estate of Harry Schutzman. D.M. is partially supported by a Brian and Tiffinie Pang Faculty Fellowship. We wish to thank Ben McKenna for an interesting discussion concerning tensor norms and for pointing us to the relevant results. We are also indebted to Mark Sellke for pointing out a few mistakes in an earlier version, as well as for suggesting a way to apply our approach to odd-degree tensors.

\section{Fourth-order tensors}\label{sec:tensors}
In the sequel, we shall mostly focus on simple examples of non-quadratic functionals in the form of fourth-order tensors. This section outlines the different perspectives we take on these tensors and presents some basic results.

\subsection{Fourth-order tensors as multi-linear forms}
Consider a fourth-order tensor \( T  \in  (\RR^n)^{\otimes 4} \), and equip $\RR^n$  with the standard Euclidean norm \( \lVert \cdot \rVert_2 := \sqrt{\langle \cdot , \cdot \rangle} \). We shall always assume that \( T \) is symmetric, i.e., 
\[
    T_{i_1 i_2 i_3 i_4} = T_{i_{\sigma(1)} i_{\sigma(2)} i_{\sigma(3)} i_{\sigma(4)}} \quad \forall \sigma \in S_4,
\]
where \( S_4 \) is the symmetric group on four elements.
Let \( \SSS^{n-1} := \{ x \in \RR^n : \lVert x \rVert_2 = 1 \} \) denote the unit sphere in \( \RR^n \). We call $(i_1,i_2,i_3,i_4)$ a \emph{diagonal} entry if it contains any repetitions, i.e. $|\{i_1,i_2,i_3,i_4\}| < 4$. The tensor $T$ can be realized as a multi-linear form $\RR^n\times\RR^n\times\RR^n\times \RR^n\to \RR$, defined on the standard basis $\{e_i\}_{i=1}^n$ through the relation
$T_{i_1 i_2 i_3 i_4} = T(e_{i_1},e_{i_2},e_{i_3},e_{i_4})$.
This perspective leads to the definition of the \emph{injective norm} of \( T \) as
\[
    \lVert T \rVert_{\inj} := \sup_{x,y,z,w \in \SSS^{n-1}} \big| T(x,y,z,w) \big|.
\]
The symmetry of $T$ affords the following, simpler, characterization of the injective norm.
\begin{Lem}{(\cite[Proposition 1]{Bochnak1971})} \label{lem:injnorm}
    Let $T\in (\RR^n)^{\otimes 4}$ be a symmetric fourth-order tensor. Then, 
    \begin{equation*}
        \lVert T \lVert_\inj = \sup_{x \in \SSS^{n-1}} \big\lvert T(x,x,x,x)\big\lvert.
    \end{equation*}
\end{Lem}
\subsection{Fourth-order tensors as matrix operators}
Above, we have treated the tensor $T$ as a multi-linear operator. There is another equivalent way to interpret fourth-order tensors; as linear operators acting on matrices. Specifically, we can ``flatten'' \( T \) into an \( n^2 \times n^2 \) matrix, i.e., an element of \( \mathrm{Mat}(\RR^{n^2}, \RR^{n^2}) \). For any matrix \( M \in \Mat(\RR^{n},\RR^{n}) \simeq \RR^{n^2}\) , the application of \( T \) yields another matrix, and we define the (matrix) operator norm of \( T \) as
\[
    \|T\|_{\op} = \sup_{\|M\|{\HS} = 1} \|TM\|{\HS} = \sup_{\|M\|{\HS} = 1} \langle M, TM \rangle{\HS}.
\]
To move between the two perspectives, we have the following identity
\[
    T(x,y,z,w) = \langle x \otimes y, T(z \otimes w) \rangle{\HS}.
\]
Above we used the Hilbert-Schmidt inner product in the space $\mathrm{Mat}(\RR^{n^2}, \RR^{n^2})$. Throughout, we shall use this same notation for the inner product, regardless of the ambient dimension, as it should be clear from the context.
With this perspective, when we say that $T$ is positive (semi-)definite, we mean that it is a positive (semi-)definite $n^2\times n^2$ matrix. Moreover, if $S$ is another symmetric $4$-tensor, we will write $T \preceq S$ to mean that $S - T$ is positive semi-definite. Similarly, $T\cdot S$ stands for the matrix product of $T$ and $S$.

Next, we establish a connection between the injective norm and the matrix-operator norm.

\begin{Lem} \label{lem:tensorbounds}
    Let \( T \) be a symmetric fourth-order tensor. Then, the following holds:
    \begin{itemize}
        \item \( \|T\|_{\inj} \leq \|T\|_{\op} \).
        \item If \( S \in \mathrm{Mat}(\RR^{n^2}, \RR^{n^2}) \) is symmetric and \( S \preceq T \) as matrices, then
        \[
            \|S\|_{\op} \leq \|T\|_{\op} \quad \text{and} \quad \|S\|_{\inj} \leq \|T\|_{\inj}.
        \]
    \end{itemize}
\end{Lem}

\begin{proof}
    For the first part, since \( \|x \otimes x\|{\HS} = 1 \) whenever \( x \in \SSS^{n-1} \), we have
    \begin{align*}
        \|T\|_{\inj} &= \sup_{x \in \SSS^{n-1}} \left| T(x,x,x,x) \right| = \sup_{x \in \SSS^{n-1}} \langle x \otimes x, T(x \otimes x) \rangle{\HS} \\
        &\leq \sup_{\|M\|{\HS} = 1} \langle M, TM \rangle{\HS} = \|T\|_{\op}.
    \end{align*}
For the second part, we know that for any matrix \( M \),
    \[
        \langle M, SM \rangle{\HS} \leq \langle M, TM \rangle{\HS},
    \]
    which implies \( \|S\|_{\op} \leq \|T\|_{\op} \). Specializing this to rank-1 matrices gives \( \|S\|_{\inj} \leq \|T\|_{\inj} \).
\end{proof}

\subsection{Concentration of the injective norm}
One important case of interest is when \( T \)  is a random Gaussian tensor. I.e., when $T$ has independent, up to symmetries, Gaussian off-diagonal entries, with all other entries being zero. In that case, the injective norm concentrates tightly around its mean, as in \cite{tensor}. For completeness, we state the next result for general degree $p$.
\begin{Thm}\label{thm:ben}
    Let $p \geq 4$ and suppose that $T$ is a symmetric tensor of order $p$ with independent, up to symmetries, $ \mathcal{N}\left(0, \frac{1}{n^{p-1}}\right) $ off-diagonal entries, and all other entries are zero. Then there exists an $ E_0(p)>0$ such that, for every $\eps > 0$
    \begin{equation*}
        \PP \left(\left|n^{\frac{p}{2}-1} \lVert T \lVert _\inj - E_0(p)\right|\leq \eps\right) \xrightarrow{n\to\infty} 0.
    \end{equation*}
    Moreover, $\frac{E_0(p)}{\sqrt{\log(p)}} \xrightarrow{p \to \infty} 1$ and $\ E_0(4) \approx 1.794$. 
\end{Thm}
\begin{proof}
 The convergence in probability is a direct consequence of \cite[Theorem 2.12]{tensor}, which concerns convergence of the ground energy. We only need to note that by symmetry, and Lemma \ref{lem:injnorm},
 $$\inf\limits_{\|x\|_2 = 1}T(x) \stackrel{\mathrm{law}}{=}-\|T\|_{\mathrm{inj}}.$$
 For the values of $E_0(p)$ see for example \cite[Lemma 2.7]{dartois2024injective}, and \cite[Remark 3.9]{dartois2024injective} for the exact value of $E_0(4)$.
\end{proof}

\def \Uvec{\mathfrak{I}}

\section{Tensorized Stochastic Localization} \label{sec:TSL}
This section introduces our variant SL process, the \emph{Tensorized Stochastic Localization} process. As explained, we focus on degree-$4$ tensors and prove Theorem \ref{thm:maindecompo}, though the new components we introduce are applicable to a wider range of functions. The proof of the decomposition theorem proceeds in several steps. First, in Section \ref{sec:firstdecomp} we show how to decompose an arbitrary symmetric $4$-tensor $T$ into $4$-tensors of the form $M\otimes M$, where $M$ is some symmetric matrix. Then, in Section \ref{sec:seconddecomp} we show how to repeat the same ideas, by successively decomposing each copy of $M$, and obtain a final rank-one decomposition, as in Theorem \ref{thm:maindecompo}.
\subsection{Decomposing \texorpdfstring{$4$}{lg}-tensors into tensor products of matrices} \label{sec:firstdecomp}
Fix $\mu$ a probability measure on $\DC$ given by $\mu(x) \propto \exp\left(T(x)\right)$, where $T \in (\RR^n)^{\otimes 4}$ is some positive definite symmetric fourth-order tensor. Fix also a function $\vphi:\DC \to \RR$. Our goal is to decompose $\mu$ into a tame mixture of simpler measures, allowing to control $\var_\mu(\vphi)$. A key component in our construction is that we allow to consider positive measures which are not necessarily normalized to be probability measures. For those measures, we recall here our definition of normalized variance, which applies to any non-negative measure $\nu$,
$$\overline{\var}_{\nu}(\vphi) = \EE_{\nu}\left[\vphi^2\right] - \frac{\EE_{\nu}\left[\vphi\right]^2 }{\EE_{\nu}\left[1\right]}.$$

The following is our first decomposition result:
\begin{Thm} \label{thm:firstdecomp}
	Let $\mu$ be a probability measure on $\DC$ given by $\mu(x) \propto \exp\left(T(x)\right)$, where $T$ is a positive definite symmetric fourth-order tensor, and let $\vphi:\DC \to \RR$. Then for every $\delta > 0$, there exists a decomposition,
	$$\mu = \int \mu_{M,R}\eta(\dd M,\dd R),$$
	where $M$ and $R$ are $n\times n$ symmetric matrices. The measures $\mu_{M,R}$ are non-negative and given by
	\begin{align*}
	\mu_{M,R}(x) &\propto \exp\left(\langle x^{\otimes 2}\otimes x^{\otimes 2}, M\otimes M\rangle\HS + \langle x,Rx\rangle \right) = \exp\left(\langle x, Mx\rangle^2 + \langle x,Rx\rangle \right).
	\end{align*}
	Moreover, we have the following properties:
	\begin{itemize}
		\item $\|M\|^2_{\op} \leq \|T\|_{\mathrm{inj}}$, $\eta$-almost surely.
		\item $\|R\|\HS \leq \delta$, $\eta$-almost surely.
		\item $\left|\var_{\mu}(\vphi)-\int\overline{\var}_{\mu_{M,R}}(\vphi)\eta\left(\dd u, \dd R\right)\right| \leq \delta$.
	\end{itemize}
\end{Thm}
Towards proving the decomposition theorem we consider the stochastic localization process, given by the following system of SDEs,
\begin{align}\label{eq:tensor_sl}
\dd F_t(x) = \langle x^{\otimes 2}-v_t,C_t\dd W_t\rangle\HS F_t(x), \quad F_0(x) = 1,
\end{align}
where $(W_t)_{t \geq 0}$ is a Dyson Brownian motion in the space $\mathrm{Mat}(\RR^n,\RR^n)$, $v_t$ is some adapted matrix-valued process to be defined soon, and $C_t$ is a $4$-tensor operating on matrices, so that $C_t\dd W_t$ is valued in $\mathrm{Mat}(\RR^n,\RR^n)$. With this SDE, we define the following measure-valued process
\begin{equation}\label{eq:def_mut}
\mu_t = F_t\mu.
\end{equation}
Before delving into the details of the process we shall first expand on the construction of the driving operator $C_t$ and the drift $v_t$. When constructing these processes we focus on the main challenges outlined in the introduction -- specifically, the term in \eqref{eq:toobig}, arising from the naïve construction, can become prohibitively large and interfere with applications. Therefore, our goal is to design $C_t$ and $v_t$ to eliminate this term without compromising the structure of the SL process.

\subsubsection{The driving matrix \texorpdfstring{$C_t$}{lg}}
Recall that $C_t$ is a linear operator, acting on matrices, and as such can be thought of as an element of $\Mat(\RR^{n^2},\RR^{n^2})$. The main role played by $C_t$ is to enforce some linear constraints along the dynamics of $\mu_t$. 
Specifically, we consider the following observables of $\mu_t$. The first ones are related to the dynamics of $\vphi$,
\begin{equation} \label{eq:phiconstraint}
\mathcal{A}_t=\EE_{\mu_t}\left[\vphi(x)(x^{\otimes 2}-v_t)\right] \text{ and }	\vphi_t= \EE_{\mu_t}\left[\vphi\right].
\end{equation}
Remark that $\vphi_t$ is a scalar and that $\mathcal{A}_t$ is matrix-valued. Similarly we define the analogous quantities for the constant function $1$,
\begin{equation} \label{eq:massconstraint}
\Uvec_t = \EE_{\mu_t}\left[x^{\otimes 2}-v_t\right] \text{ and }	\mathrm{mass}_t = \EE_{\mu_t}\left[1\right].
\end{equation}
Note that we will not enforce $\mu_t$ to remain a probability measure, and we will only guarantee that it is positive. Thus even though $\mathrm{mass}_0 = 1$, the same will not be true when $t > 0$. On the other hand, we will choose the matrix $C_t$ to facilitate a Poincar\'e inequality for positive non-normalized measures.
For this, as will soon become evident, our central quantity of interest is
\begin{equation} \label{eq:quantity}
\Phi_t := \frac{\mathrm{mass}_t\mathcal{A}_t-\vphi_t\Uvec_t}{\mathrm{mass}^{\frac{3}{2}}_t}
\end{equation}
The main role of $C_t$ will be to ensure that $\Phi_t$ remains nearly constant. We thus define the following $4$-tensor as an element of $\Mat(\RR^{n^2},\RR^{n^2})$,
\begin{equation} \label{eq:Tdef}
T_t =  T - \frac{1}{2}\int_0^tC_s^2\dd s,
\end{equation}
recalling the identification of $\Mat(\RR^{n^2},\RR^{n^2})$ with the space of $4$-tensors.
Ideally, we would like to define $C_t$ as the orthogonal projection $\mathrm{Proj}_{\mathrm{Image}(T_t)\cap \Phi_t^{\perp}}.$ 
The condition $\mathrm{Image}(C_t) \subset \mathrm{Image}(T_t)$ would ensure that that $T_t$ is a decreasing process, while having $\mathrm{Image}(C_t)$ orthogonal to $\Phi_t$, as in \cite{eldan2022spectral}, would play the same variance preservation from \eqref{eq:PI_EKZ}. The new component is that instead of keeping the expectation $\EE_{\mu_t}[\vphi]$ constant over $t$, our choice of $\Phi_t$ will instead make $\frac{\EE_{\mu_t}[\vphi]}{\mathrm{mass}_t}$ constant. As we will see, this is the correct condition that ensures that the measure $\mu_t$ satisfies a tractable Poincar\'e inequality.

At a technical level, defining $C_t$ as an orthogonal projection is highly discontinuous whenever $\Phi_t \in \mathrm{Image}(T_t)$. To countermand this issue we follow \cite{eldan2022spectral}, and define a smoothed version of the projection operator. For completeness, we provide the details of the construction in Appendix \ref{sec:projections}.
\begin{Lem}\label{fct:Cprops}(\cite[Lemma 2]{eldan2022spectral}) 
	Let $\mathcal{H}$ stand for the set of subspaces of $\RR^{n^2}$. There exists a function $C:\mathcal{H}\times \RR^{n^2} \to  \Mat(\RR^{n^2},\RR^{n^2})$ satisfying the following conditions, for any $H \in \mathcal{H}$ and $v,w \in \RR^{n^2}$:
	\begin{enumerate}
		\item $\mathrm{Image}(C(H, v))\subset H$, and when restricted to $H$, $C(H, v)$ is positive definite, and $C(H, v) \preceq \mathrm{I}_H$.
		\item The map $C(H,\cdot):\RR^{n^2}\to \Mat(\RR^{n^2},\RR^{n^2})$ is Lipschitz continuous.
		\item $\mathrm{Tr}(C(H, v)^2) \geq \mathrm{dim}(H)-1$.
		\item $\|C(H, v)v\|\leq \delta.$
	\end{enumerate}
\end{Lem}
With the above map, we now make the definition, 
\begin{equation} \label{eq:Cdef}
C_t:=C(\mathrm{Image}(T_t), \Phi_t).
\end{equation}
\subsubsection{The drift \texorpdfstring{$v_t$}{lg}}
In the previous section, we demonstrated how the definition of $C_t$ can be adapted to incorporate the time-varying quantity $\mathrm{mass}_t$. The total mass $\mathrm{mass}_t$ becomes relevant since we no longer consider the barycenteric process $\aa_t$ from \eqref{eq:NaiveSL}, which ensures that $\mu_t$ is a probability measure. One desirable effect of removing $\aa_t$ is that the process $\mathtt{L}_t$ from \eqref{eq:toobig} immediately becomes smaller. Recall that allowing $\mathtt{L}_t$ to attain large values was a significant hurdle in generalizing the stochastic localization process. To further reduce $\mathtt{L}_t$, even to arbitrarily small values, we capitalize on the newly gained flexibility that arises from allowing $\mathrm{mass}_t$ to change. Our next key idea is to replace the process $\aa_t$ with another adapted process that ensures $\mathtt{L}_t$ remains as small as desired.

We prove the following Proposition in Appendix \ref{sec:boundedproc}. The construction relies on ideas stemming from Bessel processes, which are known to remain positive almost surely. With this result, we define $v_t := v_t^\delta.$
\begin{Prop} \label{prop:smallv}
	Let $C_t$ be defined by $\eqref{eq:Cdef}$, and let $\delta > 0$. There exists an adapted drift process $v_t^\delta \in \mathrm{Mat}(\RR^{n},\RR^n)$, such that if $X^{\delta}_t$ satisfies
	\begin{equation}
	\label{eq:boundedprocdef1}
	\dd X^{\delta}_t = C_t\dd W_t - C^2_tv^{\delta}_t\dd t, \ \ \ X^\delta_0 = 0,
	\end{equation}
	then $\sup\limits_{t \geq 0}\|X^{\delta}_t\|\HS \leq \delta$ almost surely.
\end{Prop}
\subsubsection{Comparison with the standard SL process}
Before proceeding with the proof, we highlight the main modifications made to the original stochastic localization process. Let $\nu_0$ be a probability measure on $\DC$ and let $(B_t)_{t \geq 0}$ be a standard Brownian motion on $\mathbb{R}^n$ with $B_0 = 0$. Recall the following system of SDEs:
\begin{equation*}
\dd \FF_t(x) = \langle x - \aa_t, \CC_t \dd B_t \rangle, \quad \FF_0(x) = 1,
\end{equation*}
where standard SL is defined as the measure-valued process $\nu_t = \FF_t \nu_0$. Here, recall $\aa_t$ is the barycenter,
\[
\aa_t = \int_{\DC} x \, \dd \nu_t(x),
\]
and $\CC_t$ is the driving matrix. Let us highlight the main differences between the two processes.
\begin{enumerate}
	\item Tensorized Stochastic Localization (TSL) is defined for quadratic forms $x ^{\otimes 2}$ rather than on $x$, as in SL. This is adapted to have a measure with higher-order potentials. In fact, let us define for a matrix $A_t$ and a vector $b_t$, both depending on $t$, the operators
	\begin{equation*}
	Q(A_t):= \int_0^t A_s^2 \dd s \quad \text{ and } \quad L(A_t, b_t, B_t) = \int_0^t A_s \dd B_s-\int_0^t A_s^2 b_s \dd s.
	\end{equation*}
	Then, for SL $\FF_t$ has the form \cite[Fact 14]{ChenEldan}
	\begin{equation*}
	\FF_t(x) \propto \exp \left(-\frac{1}{2}\langle Q(\CC_t)x,x\rangle+\langle x, L(\CC_t, \aa_t, B_t)\rangle \right)\,,
	\end{equation*}
	instead for TSL, see Lemma \ref{lem:densityexp}, $F_t$ has the form 
	\begin{equation*}
	F_t(x) \propto \exp\left(T(x)-\frac{1}{2}\langle Q(C_t)x^{\otimes 2},x^{\otimes 2}\rangle\HS+\langle x^{\otimes 2}, L(C_t, v_t, W_t) \rangle\HS \right)\,.
	\end{equation*}
	\item In SL the driving process $B_t$ is a standard Brownian motion. Instead in TSL the driving process $W_t$ is a Dyson Brownian motion, which evolves in the space of symmetric matrices. This is essential because more randomness in the quadratic variation is required to generate a degree-$4$ tensor for the decomposition, alongside a quadratic term.
	\item The quadratic term is canceled out by the drift $v_t$. In contrast, in SL the role of the barycenter $\aa_t$ is to ensure that $\nu_t$ remains a probability measure. As a result, $\mu_t$ is no longer a probability measure, but an unnormalized positive measure. In TSL, this is reflected by incorporating a new linear constraint to the driving matrix $C_t$. 
	\item Furthermore, since $W_t$ is an $n \times n$ matrix, $C_t$ now acts as an operator on matrices, which can also be interpreted as a fourth-order tensor.
\end{enumerate}

We also mention that there is another possible path to take. Instead of defining $C_t$ in a way that will later enable us to handle un-normalized measures, we could encode a linear constraint that makes $\mathrm{mass}_t = \EE_{\mu_t}[1]$ constant over time. In that case, we would ensure that $\mu_t$ is a probability measure, with the role of the barycenter $\aa_t$ now incorporated into $C_t$. In fact, this is the choice we made in a previous version of this work. However, adding an extra constraint into $C_t$ resulted in a decomposition into low-rank tensors with rank strictly higher than $1$. This fact has led to a deterioration in the obtained spectral gaps. We do believe that this approach could be of interest in other settings though. 

\subsubsection{Properties of \texorpdfstring{$\mu_t$}{lg}} With the above definitions at hand, we are ready to prove some key properties of the measures $\mu_t$. We first observe that while $\mu_t$ does not remain a probability measure, its normalized variance $\overline \var_{\mu_t}(\vphi)$ does not vary too much through the process.
\begin{Lem}\label{lem:new}
	Let $(\mu_t)_t$ be defined as in \eqref{eq:def_mut}. Then, for every $t \geq 0$,
		$$\left|\EE\left[\overline{\var}_{\mu_t}(\vphi)\right] - \var_\mu(\vphi)\right| \leq t\delta^2.$$
\end{Lem}
\begin{proof}
	We begin by representing $\frac{\vphi^2_t}{\mathrm{mass}_t}$ as an It\^o process.
	By definition of $\mu_t$, in \eqref{eq:tensor_sl}, 
	\begin{align} \label{eq:volcalc}
	\dd\vphi_t &= \dd\int\vphi(x) F_t(x)\mu(\dd x) = \int \vphi(x)\langle x^{\otimes 2}-v_t,C_t\dd W_t\rangle\HS F_t(x)\mu(\dd x) \nonumber\\
	&= \left\langle \int \vphi(x)(x^{\otimes 2}-v_t )F_t(x)\mu(\dd x),C_t\dd W_t\right\rangle\HS = \langle C_t\mathcal{A}_t,\dd W_t\rangle\HS,\nonumber\\
	\dd\mathrm{mass}_t &= \dd\int F_t(x)\mu(\dd x) = \int \dd F_t(x)\mu(\dd x) = \int \langle x^{\otimes 2}-v_t,C_t\dd W_t\rangle\HS F_t(x)\mu(\dd x) \nonumber\\
	&= \left\langle \int (x^{\otimes 2}-v_t )F_t(x)\mu(\dd x),C_t\dd W_t\right\rangle\HS = \langle C_t\Uvec_t,\dd W_t\rangle\HS,
	\end{align}
	where $\mathcal{A}_t$ and $\Uvec_t$ are defined in \eqref{eq:phiconstraint} and \eqref{eq:massconstraint}. Applying It\^o's lemma to the function $(x,y)\to \frac{x^2}{y}$, we get
	\begin{align*}
	\dd \frac{\vphi^2_t}{\mathrm{mass}_t} &= 2\frac{\vphi_t}{\mathrm{mass}_t}\dd\vphi_t + \frac{1}{\mathrm{mass}_t}\dd[\vphi]_t - \frac{\vphi^2_t}{\mathrm{mass}^2_t}\dd\mathrm{mass}_t + \frac{\vphi^2_t}{\mathrm{mass}^3_t}\dd[\mathrm{mass}]_t - 2\frac{\vphi_t}{\mathrm{mass}_t^2}d[\mathrm{mass},\vphi]_t\\
	&=2\frac{\vphi_t}{\mathrm{mass}_t}\dd\vphi_t  - \frac{\vphi^2_t}{\mathrm{mass}^2_t}\dd\mathrm{mass}_t + \frac{1}{\mathrm{mass}_t}\|C_t\mathcal{A}_t\|\HS^2 - 2\frac{\vphi_t}{\mathrm{mass}^2_t}\langle C_t\Uvec,C_t\mathcal{A}_t\rangle\HS + \frac{\vphi^2}{\mathrm{mass}_t^3}\|C_t\Uvec_t\|\HS^2\\
	&= 2\frac{\vphi_t}{\mathrm{mass}_t}\dd\vphi_t  - \frac{\vphi^2_t}{\mathrm{mass}^2_t}\dd\mathrm{mass}_t + \left\|C_t\left(\frac{\mathrm{mass}_t\mathcal{A}_t-\vphi_t\Uvec_t}{\mathrm{mass}_t^{\frac{3}{2}}}\right)\right\|\HS^2\\
	&= 2\frac{\vphi_t}{\mathrm{mass}_t}\dd\vphi_t  - \frac{\vphi^2_t}{\mathrm{mass}^2_t}\dd\mathrm{mass}_t + \left\|C_t\Phi_t\right\|\HS^2\\
	& \leq 2\frac{\vphi_t}{\mathrm{mass}_t}\dd\vphi_t  - \frac{\vphi^2_t}{\mathrm{mass}^2_t}\dd\mathrm{mass}_t + \delta^2,
	\end{align*}
	where we have used the expression \eqref{eq:phiconstraint} for $\Phi_t$, the definition of $C_t$ from \eqref{eq:Cdef} and Lemma \ref{fct:Cprops}.
	Observe that from \eqref{eq:volcalc} both $\vphi_t$ and $\mathrm{mass}_t$ are martingales. Thus, since $\mu_0 = \mu$ and $\mathrm{mass}_0 = 1$,
	$$\left|\EE\left[\frac{\vphi^2_t}{\mathrm{mass}_t}\right]- \EE\left[\frac{\vphi^2_0}{\mathrm{mass}_0}\right]\right| = \left|\EE\left[\frac{\EE_{\mu_t}\left[\vphi\right]^2}{\EE_{\mu_t}\left[1\right]}\right] - \EE_{\mu}\left[\vphi\right]^2\right| =\left|\int_0^t\frac{\dd}{\dd s}\EE\left[\frac{\vphi^2_t}{\mathrm{mass}_t}\right]\dd s\right|\leq t\delta^2.$$
	Furthermore $\EE_{\mu_t}\left[\vphi^2\right]$ is also a martingale and we conclude with 
	\begin{align*}
	\left|\EE\left[\overline{\var}_{\mu_t}(\vphi)\right] - \var_\mu(\vphi)\right| = \left|\EE\left[\EE_{\mu_t}\left[\vphi^2\right]\right]-\EE_{\mu}\left[\vphi^2\right] + \EE\left[\frac{\EE_{\mu_t}\left[\vphi\right]^2}{\EE_{\mu_t}\left[1\right]}\right] - \EE_{\mu}\left[\vphi\right]^2\right| \leq t\delta^2.
	\end{align*}
\end{proof}
We next show that $\mu_t$ remains a positive measure and determine its (un-normalized) density.
\begin{Lem} \label{lem:densityexp}
	Almost surely, for every $t \geq 0$, $\mu_t$ is a non-negative measure on $\DC$, with the following density,
	\begin{equation} \label{eq:mudensity}
	\mu_t(x) \propto \exp\left(T_t(x) + \left\langle x^{\otimes 2}, \int\limits_0^tC_s\dd W_s -\int\limits_0^tC_s^2v_s\dd s\right\rangle\HS\right),
	\end{equation}
	where $T_t$ is given by \eqref{eq:Tdef}.
\end{Lem}
\begin{proof}
	By It\^o's lemma, we have
	\begin{equation*}
	\dd \log F_t(x) = \langle x^{\otimes 2}-v_t, C_t\dd W_t\rangle\HS -\frac{1}{2}\lVert C_t(x^{\otimes 2}-v_t) \lVert^2\HS \dd t\,.
	\end{equation*} 
	Integrating we get 
	\begin{equation*}
	F_t(x) \propto \exp \left( \left\langle x ^{\otimes 2}, \int_{0} ^{t} C_s \dd W_s - \int_{0}^{t} C_s^2 v_s \dd s \right\rangle\HS - \frac{1}{2}\left\langle x^{\otimes 2},\left(\int_{0}^{t} C_s^2 \dd s\right)x^{\otimes 2}\right\rangle\HS\right)
	\end{equation*} and by \eqref{eq:Tdef} we conclude.
\end{proof}

Finally, before going into the proof of the decomposition, we show that the process $T_t$, defined in \eqref{eq:Tdef} is decreasing.
\begin{Lem} \label{lem:Tprops}
	Let $T_t$ be the process defined in \eqref{eq:Tdef}. Then $T_t$ is a decreasing process of positive semi-definite matrices in $\Mat(\RR^{n^2},\RR^{n^2})$.
	Moreover, if $\tau = \inf\{t\geq 0| \mathrm{rank}(T_t) \leq 1\}$ then almost surely $\tau \leq 2\mathrm{Tr}(T).$
\end{Lem}
\begin{proof}
	By Lemma \ref{fct:Cprops} it holds that $\int_0^tC_s^2\dd s$ is positive semi-definite. Hence, since  $T_t = T - \frac{1}{2}\int_0^tC_s^2\dd s$ it is clear that $T_t$ is decreasing. Moreover,
	$\frac{\dd T_t}{\dd t} = -\frac{1}{2} C_t^2$, and by construction in \eqref{eq:Cdef}, and Lemma \ref{fct:Cprops}, $\mathrm{Image}(C^2_t) \subset \mathrm{Image}(T_t)$. Thus, if any $v \in \RR^{n^2}$, satisfies for some $t' \geq 0$ that $T_{t'}v = 0$, then $T_{t}v= 0$, for every $t \geq t'$ as well. Since $T_0 = T$ is positive semi-definite, this property is maintained throughout the process.
	Finally, notice that, for any $t < \tau$, again by Lemma \ref{fct:Cprops},
	$$\frac{\dd}{\dd t}\mathrm{Tr}(T_t) = -\frac{1}{2}\mathrm{Tr}(C^2_t) \leq \frac{1- \mathrm{dim}(\mathrm{Image}(T_t))}{2}\leq -\frac{1}{2}.$$
	Thus, suppose that $\tau > \mathrm{Tr}(T)$, we would then get that $\mathrm{Tr}(T_t) < 0$, which cannot happen when $T_t$ is positive semi-definite.
\end{proof}
\subsubsection{Proof of the decomposition theorem}
\begin{proof}[Proof of Theorem \ref{thm:firstdecomp}]
	Fix $\delta > 0$, and consider the stopping time $\tau = \inf\{t\geq 0| \mathrm{rank}(T_t) \leq 1\}$, as in Lemma \ref{lem:Tprops}.
	Since $\mu_t$ is a martingale, we have that $\mu = \EE\left[\mu_\tau\right]$, where the expectation is taken over the randomness of the process.
	Since $\mathrm{rank}(T_\tau)\leq 1$, we can rewrite it as $T_\tau = M\otimes M$, where $M$ is a symmetric $n \times n$ matrix. By Lemma \ref{lem:Tprops} $T_\tau \preceq T$, and we conclude that, as in Lemma \ref{lem:tensorbounds},
	$$\|M\|_{\op}^2 = \|M\otimes M\|_{\mathrm{inj}}\leq \|T_\tau\|_{\mathrm{inj}}\leq \|T\|_{\mathrm{inj}}.$$
	We also write $R = \int\limits_0^\tau C_t\dd W_t -\int\limits_0^\tau C^2_tv_t\dd t$, and observe that by Proposition \ref{prop:smallv} and the choice of $v_t$, $\|R\|\HS \leq \delta.$ Thus, by Lemma \ref{lem:densityexp}
	\begin{align} \label{eq:approxdecomp}
	\mu(x) = \EE\left[\mu_\tau(x)\right] \propto \EE\left[\exp\left(\langle x^{\otimes 2}\otimes x^{\otimes 2}, M\otimes M\rangle\HS + \langle x^{\otimes 2}, R \rangle\HS \right)\right],
	\end{align}
	which is the desired decomposition, with the mixing measure as the implicit joint distribution of $(M,R)$.
	As for the variance of $\vphi$, by Lemma \ref{lem:new} and Lemma \ref{lem:Tprops}, we get
	$$\left|\EE\left[\overline{\var}_{\mu_\tau}(\varphi)\right] -\var_{\mu}(\varphi)\right| \leq \EE\left[\tau\right]\delta^2 \leq 2\mathrm{Tr}(T)\delta^2 .$$
	We can obtain the final bound by replacing $\delta$ with $\frac{\delta}{2\mathrm{Tr}(T)+1}$.
\end{proof}

\subsection{Decomposing matrix products} \label{sec:seconddecomp}
In light of Theorem \ref{thm:firstdecomp}, our next goal will be to better understand measures $\nu$ of the form,
$$\nu(x) \propto \exp\left(\langle x,Mx\rangle^2\right),$$
where $M$ is a symmetric matrix. 
The following decomposition theorem will be our main tool for studying such measures.
\begin{Thm} \label{thm:seconddecomp}
	Let $\tilde{\nu}$ be a reference non-negative measure on $\DC$ and let $M',M$ be two positive definite matrices. Define the non-negative measure $\nu$ by
	$$\nu(x) = \exp(\langle x,M'x\rangle\langle x,Mx\rangle)\tilde{\nu}(x).$$
	Then, for every $\delta> 0$, there exists a decomposition,
	$$\nu = \int \nu_{u,R}\eta(\dd u,\dd R),$$
	where $u,R \in \RR^n$ and $\nu_{u,R}$ are non-negative measures given by,
	$$\nu_{u,R}(x) \propto \exp\left(\langle x,M'x\rangle \langle x, u\rangle^2+ \sqrt{\langle x, M'x\rangle}\langle x,R\rangle
	\right)\tilde{\nu}(x).$$
	Moreover, we have the following properties:
	\begin{itemize}
		\item $\|u\|^2 \leq \|M\|_{\op}$ $\eta$-almost surely.
		\item $\|R\| \leq \delta$ $\eta$-almost surely.
		\item $\left|\overline{\var}_{\nu}(\vphi)-\int\overline{\var}_{\nu_{u,R}}(\vphi)\eta\left(\dd u, \dd R\right)\right| \leq \delta$.
	\end{itemize}
\end{Thm}
We take a similar approach to the proof of Theorem \ref{thm:firstdecomp}, and define a variant of the stochastic localization process, taking into account the structure of $\nu$.

Thus define the following process on relative densities (we use the overlines to emphasize the overline to distinguish from the analog processes which appear in the proof of Theorem \ref{thm:firstdecomp}), 
$$\dd F_t(x) = \sqrt{\langle x,M'x\rangle} \langle x-\overline v_t,\overline C_t\dd B_t\rangle F_t(x), \quad F_0(x) = 1,$$
and set
$$\nu_t = F_t\nu.$$
Note that $B_t$ is a standard Brownian motion in $\RR^n$, unlike the Dyson Brownian motion in \eqref{eq:tensor_sl}. The driving matrix $\overline C_t$, now an element in $\mathrm{Mat}(\RR^n,\RR^n)$ is similar to the one constructed in \eqref{eq:Cdef}. We define two vector-valued process
$$\overline{\mathcal{A}}_t=\EE_{\nu_t}\left[\vphi(x)\sqrt{\langle x,M'x\rangle}(x-v_t)\right], \ \ \ \overline\Uvec_t = \EE_{\nu_t}\left[\sqrt{\langle x,M'x\rangle}(x-v_t)\right],$$
with their corresponding expectation processes
$$\overline{\vphi}_t=\EE_{\nu_t}\left[\vphi\right], \ \ \ \overline{\mathrm{mass}}_t = \EE_{\nu_t}\left[1\right],$$
and the target vector
$$\overline{\Phi}_t = \frac{\overline{\mathrm{mass}_t}\overline{\mathcal{A}}_t-\overline{\vphi}_t\overline{\Uvec}_t}{\overline{\mathrm{mass}_t}^{\frac{3}{2}}}.$$
As before we will consider the decreasing matrix-valued potential process,
\begin{equation} \label{eq:Mt}
M_t := M - \frac{1}{2}\int\limits_0^t \overline C_s^2\dd s.
\end{equation}
The matrix $\overline C_t$ is the smoothed version of $\mathrm{Proj}_{\mathrm{Image}(M_t)\cap \overline\Phi_t^{\perp}}$,
$$\overline C_t = C(\mathrm{Image}(M_t),\overline\Phi_t),$$
where the mapping $C$ is defined as in \eqref{eq:Cdef}, with the obvious modifications. In particular, Lemma \ref{fct:Cprops} applies, \emph{mutatis mutandis}.
The vector-valued drift $\overline v_t$ again serves the same role as in Theorem \ref{thm:firstdecomp} and we choose it so that the process defined by \begin{equation}\label{eq:proddrift}
\dd X_t = \overline C_t\dd B_t - \overline C_t^2\overline v_t\dd t, \ \ X_0=0,
\end{equation}
satisfies $\|X_t\| \leq \delta,$ almost surely for every $t \geq 0$. The existence of $\overline v_t$ is proved in Section \ref{sec:boundedproc}.

Again we show that the normalized variance  $\overline{\var}_{\nu_t}(\varphi)$ is a conserved quantity.
\begin{Lem} \label{lem:proddynamics}
	Let $(\nu_t)_t$ be defined as above. Then, for every $t \geq 0$,
	$$		\left|\EE\left[\overline{\var}_{\nu_t}(\vphi)\right] - \overline\var_\nu(\vphi)\right| \leq t\delta^2,$$
\end{Lem}

\begin{proof}
	The proof is essentially identical to the proof of Lemma \ref{lem:new}.
	Since $\overline{\vphi}_t$ and $\overline{\mathrm{mass}}_t$ are martingales, the same calculations show that
	$$\frac{\dd}{\dd s} \EE\left[\frac{\overline{\vphi}^2_t}{\overline{\mathrm{mass}}_t}\right] \leq \delta.$$
	The final bound then follows from the same argument.
\end{proof}

We also obtain a similar formula for the density of $\nu_t$.
\begin{Lem} \label{lem:proddensityexp}
	Almost surely, for every $t \geq 0$, $\nu_t$ is a non-negative measure on $\DC$, with the following density,
	\begin{equation}
	\nu_t(x) \propto \exp\left(\langle x,M'x\rangle\langle x,M_tx\rangle + \sqrt{\langle x,M'x\rangle}\left\langle x, \int\limits_0^t\overline C_s\dd W_s -\int\limits_0^t\overline C_s^2\overline v_s\dd s\right\rangle\right)\tilde{\nu}(x),
	\end{equation} 
	where $M_t$ is given by \eqref{eq:Mt}.
\end{Lem}
\begin{proof} Similar to the proof of Lemma \ref{lem:densityexp}, we apply It\^o's lemma and we get
	\begin{align*}
	\dd \log F_t(x) &= \frac{\dd F_t(x)}{F_t(x)}-\frac{1}{2}\frac{\dd [F_t(x)]_t}{F_t(x)^2}\\
	&= \sqrt{\langle x, M'x\rangle}\langle x-v_t, C_t\dd B_t\rangle -\frac{1}{2}\langle x, M'x\rangle \lVert C_t(x-v_t)\lVert^2 \dd t.
	\end{align*}
	By integrating and \eqref{eq:Mt} we conclude. 
\end{proof}
Let us now prove Theorem \ref{thm:seconddecomp}
\begin{proof}[Proof of Theorem \ref{thm:seconddecomp}]
	We follow the same steps as in the proof of Theorem \ref{thm:firstdecomp}. First observe that by definition of $\overline C_t$, the process $M_t$ in \eqref{eq:Mt} is a decreasing process of positive semi-definite matrices.
	We thus define the stopping time $\tau = \inf\{t\geq 0|\mathrm{rank}(M_t) \leq 1\}$ and observe that for any $t < \tau$, $\mathrm{Tr}(C_t^2) \geq 1$. Since $\frac{\dd}{\dd t}\mathrm{Tr}(M_t) = -\frac{1}{2}\mathrm{Tr}(C_t^2)$ and since $M_t$ is positive semi-definite, we can conclude that $\tau \leq 2\mathrm{Tr}(M),$ almost surely.
	Necessarily, $\mathrm{rank}(M_\tau)\leq 1$, and we can write, $M_\tau = U\otimes U,$ for some random vector $U \in \RR^n$.
	Also, because $M_t$ is decreasing, we get that,
	$$\|U\|^2\leq \|M\|_{\op}.$$
	By the optional stopping theorem, and the expression in Lemma \ref{lem:proddensityexp}, we obtain the decomposition,
	$$\nu(x) = \EE\left[\nu_\tau(x)\right] = \EE\left[\exp\left(\langle x, M'x\rangle\langle U,x\rangle^2  +\sqrt{\langle x, M'x\rangle}\langle x,R\rangle \right)\tilde{\nu}(x)\right],$$
	where $R = \int\limits_0^t\overline C_t\dd B_t -\int\limits_0^t\overline C^2_t\overline v_t\dd t,$ and satisfies 
	$\|R\| \leq \delta,$ by the choice in \eqref{eq:proddrift}.
	$$		\left|\EE\left[\overline{\var}_{\nu_t}(\vphi)\right] - \overline\var_\nu(\vphi)\right| \leq t\delta^2.$$
    The proof concludes by appropriately substituting $\delta = \delta\min\left(\frac{1}{\sqrt{n\|M\|_{\mathrm{op}}}+1},\frac{1}{2\mathrm{Tr}(M)+1}\right).$
\end{proof}
\subsection{A decomposition into low-rank tensors - Proof of Theorem \ref{thm:maindecompo}}
With some repeated applications of Theorem \ref{thm:seconddecomp}, we can now handle measures obtained through the decomposition in Theorem \ref{thm:firstdecomp}, and thus prove Theorem \ref{thm:maindecompo}. The only remaining component is that the matrices produced by Theorem \ref{thm:firstdecomp} are not positive definite, a necessary condition to apply Theorem \ref{thm:seconddecomp}. This fact is easily handled by shifting the matrices, and applying the decomposition theorem of \cite{eldan2022spectral}, i.e. Theorem \ref{thm:rank1decompo}.
\begin{proof}[Proof of Theorem \ref{thm:maindecompo}]
	We begin by applying Theorem \ref{thm:firstdecomp} for a decomposition,
	\begin{equation} \label{eq:firstdecomp}
	\mu = \int \mu_{M,R}\eta'(\dd M,\dd R),
	\end{equation}
	where 
	$$\mu_{M,R}(x) \propto \exp\left(\langle x^{\otimes 2}\otimes x^{\otimes 2}, M\otimes M\rangle + \langle x, Rx\rangle\right) = \exp\left(\langle x,Mx\rangle ^2 + \langle x,Rx\rangle\right),$$
	with $\|M\|^2_\op\leq \|T\|_{\mathrm{inj}}, \|R\|\HS\leq \delta$, $\eta'$-almost surely.
	$M$ is not necessarily positive definite, but the above shows that,
	$$0 \preceq \tilde M := M + \sqrt{\|T\|_{\mathrm{inj}}}\mathrm{I}_n\preceq 2\sqrt{\|T\|_{\mathrm{inj}}}\mathrm{I}_n,$$ is almost surely positive semi-definite. To rewrite the density in a useful way, we observe
	\begin{align*}
	\langle x,Mx\rangle^2 &= \left\langle x,\left(\tilde M - \sqrt{\|T\|_{\mathrm{inj}}}\mathrm{I}_n\right)x\right\rangle^2\\
	&= \langle x,\tilde Mx\rangle ^2 + \|T\|_{\mathrm{inj}}\|x\|^4 -2\sqrt{\|T\|_{\mathrm{inj}}}\|x\|^2\langle x,\tilde Mx\rangle\\
	&= \langle x,\tilde Mx\rangle ^2 + \|T\|_{\mathrm{inj}}n^2 -2\sqrt{\|T\|_{\mathrm{inj}}}n\langle x,\tilde Mx\rangle.
	\end{align*}
	Thus, we obtain the following identity
	$$ \mu_{M,R}(x) \propto \exp\left(\langle x,\tilde{M}x\rangle ^2 + \langle x,Rx\rangle-2n\sqrt{\|T\|_{\mathrm{inj}}}\langle x,Mx\rangle\right).$$
	To apply Theorem \ref{thm:seconddecomp}, we isolate $\tilde{M}$ and define,
	$$\tilde{\mu}_{M,R} \propto \exp\left(\langle x, Rx\rangle-2n\sqrt{\|T\|_{\mathrm{inj}}}\langle x,\tilde{M}x\rangle\right)= \exp\left(\langle x, Rx\rangle\right)\bar{\mu}_{M}(x),$$
	where 
	\begin{equation} \label{eq:mbar}
	\bar{\mu}_{M}(x) \propto \exp\left(-2n\sqrt{\|T\|_{\mathrm{inj}}}\langle x,\tilde{M}x\rangle\right).
	\end{equation}
	With these definitions,
	$$\mu_{M,R}(x) \propto \exp\left(\langle x, \tilde{M}x\rangle\langle x, \tilde{M}x\rangle\right)\tilde{\mu}_{M,R}(x).$$
	Now, by Theorem \ref{thm:seconddecomp}, applied to a single copy of $\tilde{M}$, we obtain the decomposition,
	$$\mu_{M,R}(x) \propto \int\exp \left(\langle x,\tilde{M}x\rangle\langle u,x\rangle^2 + \sqrt{\langle x, \tilde{M}x\rangle}\langle x,R'\rangle\right)\tilde{\mu}_{M,R}(x)\eta_1(\dd u,\dd R'),$$
	where $u, R' \in \RR^n$ satisfy $\|u\|^2 \leq \|\tilde{M}_1\|_\op\leq 2\sqrt{\lVert T\lVert_\inj}$, and $\|R\|\leq \delta$, $\eta_1$-almost surely.
We then set $\mu_{M,u,R,R'} = \exp \left(\langle x,\tilde{M}x\rangle\langle u,x\rangle^2 + \sqrt{\langle x, \tilde{M}x\rangle}\langle x,R'\rangle\right)\tilde{\mu}_{M,R}(x)$ and note that (with the understanding that $\eta_1$ depends on $M$ and $R$),
	\begin{align}  \label{eq:iteratingdelta}
	&\left|\var_{\mu}(\vphi)-\int\overline{\var}_{\mu_{M,u,R,R'}}(\vphi)\eta'(\dd M,\dd R)\eta_1(\dd u,\dd R')\right| \nonumber\\
	&\ \ \ \leq \left|\var_{\mu}(\vphi)-\int\overline{\var}_{\mu_{M,u}}(\vphi)\eta'(\dd M,\dd R)\right|\nonumber\\
	&\ \ \ \ \ + \left|\int\left(\overline{\var}_{\mu_{M,u}}(\vphi)-\int\overline{\var}_{\mu_{M,u,R,R'}}(\vphi)\eta_1(\dd u,\dd R')\right)\eta'(\dd M,\dd R) \right|\nonumber\\
	&\ \ \ \leq 2\delta.
	\end{align}
	We proceed by writing 
	$$\langle u,x\rangle^2 = \langle x, (u \otimes u) x\rangle,$$
	and applying Theorem \ref{thm:seconddecomp} again on each $\mu_{M,u,R,R'}$,
	separately to obtain a refined decomposition
	$$\mu_{M,R}(x) \propto \int\exp \left(\langle u,x\rangle^2\langle v,x\rangle^2 + \sqrt{\langle x,\tilde{M}x\rangle}\langle x,R'\rangle + \sqrt{\langle u,x\rangle^2}\langle x,R''\rangle \right)\tilde{\mu}_{M,R}(x)\eta_2(\dd u,\dd v,\dd R', \dd R''),$$
	with $\|u\|^2,\|v\|^2 \leq \|\tilde{M}_1\|_\mathrm{op} \leq 2\sqrt{\lVert T\lVert_\inj}$ and  $\|R'\|,\|R''\|\leq \delta$, $\eta_2$-almost surely.

	To get rid of the quadratic part we switch roles with $\bar{\mu}_{M}$, from \eqref{eq:mbar}, and use Theorem \ref{thm:rank1decompo} to get,
	\begin{align} \label{eq:finaldecomp}
	\mu(x) 
	=& \int\exp\left(\langle u,x\rangle^2\langle v,x\rangle^2+\langle w, x\rangle^2 + \langle \ell, x\rangle + \psi(x)\right)\eta(\dd u, \dd v, \dd w,\dd\ell,\dd \psi),
	\end{align}
	where $\psi(x) =  \sqrt{\langle x,\tilde{M}x\rangle}\langle x,R'\rangle + \sqrt{\langle u,x\rangle^2}\langle x,R''\rangle + \langle x,Rx\rangle^2$ satisfies, by replacing $\delta$ with an appropriate choice, $\max\limits_{x\in \DC}|\psi(x)| \leq\delta$ and where $\|w\|^2\leq 2n\sqrt{\|T\|_{\mathrm{inj}}}\cdot\|\tilde{M}\|_{\mathrm{op}}\leq 4n\|T\|_{\mathrm{inj}}$. The decomposition of the variance follows straightforwardly by iterating the argument in \eqref{eq:iteratingdelta}.
\end{proof}
\section{Spectral gap estimates} \label{sec:spectral_gap}
We now apply our Theorem \ref{thm:maindecompo} to establish spectral gaps for tensor Ising models. We will primarily use our decomposition result to reduce the problem from general tensor models to rank-one tensors and then prove Theorem \ref{thm:tensorgap}. Thus, we first need to collect some facts about the spectral gaps of these simpler measures.
\subsection{Spectral gaps for simple measures - Dobrushin's Condition}
Our goal in this section is to study the spectral gap of simple measures that are inherently low rank. The primary tool we will use to analyze these quantities is the \emph{Dobrushin condition} \cite{Dobruschin1968}. To state this condition, recall that if $\mu$ is a measure on $\DC$ and $X \sim \mu$, then the influence matrix of $\mu$ is defined as
\begin{equation*}
A^\mu_{ij} = \frac{1}{2} \sup_{x^i, \tilde{x}^i} \left| \EE\left[X_i | x^i\right] - \EE\left[X_i | \tilde{x}^i\right] \right|,
\end{equation*}
where the supremum is taken over all pairs $x^i, \tilde{x}^i \in \{-1, 1\}^{[n] \setminus \{i\}}$ such that $x^i$ and $\tilde{x}^i$ differ only in the $j$-th bit. The Dobrushin condition holds when $\|A^{\mu}\|_{{\op}} < 1$ and implies that
\begin{equation} \label{eq:dobru}
C_{\poi}(\mu) \leq \frac{1}{1 - \|A^{\mu}\|_{{\op}}}.
\end{equation}
To facilitate bounds on the influence matrix, we introduce an auxiliary construct, which we call the \emph{derivative matrix}. For a function $F: \DC \to \mathbb{R}$, the derivative matrix is defined as
\begin{equation*}
D^F_{ij} = \frac{1}{4} \sup_{x^i, \tilde{x}^i} \left| F(x^i_+) - F(x^i_-) - \left( F(\tilde{x}^i_+) - F(\tilde{x}^i_-) \right) \right|,
\end{equation*}
where $x^i_\pm$ denotes the vector $x^i$ with the $i$-th bit set to $\pm 1$. To see the connection to derivatives, consider the discrete derivative of $F$ on the hypercube
\begin{equation*}
\partial_i F(x) =\frac{F(x^i_+) - F(x^i_-)}{2}.
\end{equation*}
With this definition
\begin{equation*}
D^F_{ij} = \max_{x} |\partial_j \partial_i F(x)|.
\end{equation*}

We now show that the derivative matrix provides a bound on the influence matrix and that we can establish these bounds in an additive manner.
\begin{Lem} \label{lem:additivedob}
Let $\{F_k\}_{k=1}^K: \DC \to \mathbb{R}$ and define

\begin{equation*}
\mu(x) \propto \exp\left(\sum_{k=1}^K F_k(x)\right).
\end{equation*}
Then,
\begin{equation*}
\|A^{\mu}\|_{{\op}} \leq \sum_{k=1}^K \|D^{F_k}\|_{{\op}}.
\end{equation*}
\end{Lem}
\begin{proof}
Let $i,j \in [n]$ and note that if $X \sim \mu$, for any $x^i \in \{-1,1\}^{[n] \setminus\{i\}}$, we have that $X_i|x^i$ is a Bernoulli random variable with parameter $p$, where $\frac{p}{1-p} = \exp\left\{\sum_{k=1}^K \left(F_k(x_+^i) - F_k(x_-^i)\right)\right\}$. A quick calculation shows that
\begin{equation*}
\EE\left[X_i | x^i\right] = \tanh\left(\frac{1}{2}\log\left(\frac{p}{1-p}\right)\right).
\end{equation*}
Since $\tanh$ is 1-Lipschitz, we obtain
\begin{align*}
A^{\mu}_{ij} &= \frac{1}{2}\sup_{x^i,\tilde{x}^i} \left|\EE\left[X_i|x^i\right] - \EE\left[X_i|\tilde{x}^i\right]\right| \\
&= \frac{1}{2}\sup_{x^i,\tilde{x}^i}\left|\tanh\left(\frac{1}{2}\sum_{k=1}^K \left(F_k(x_+^i) - F_k(x_-^i)\right)\right) - \tanh\left(\frac{1}{2}\sum_{k=1}^K \left(F_k(\tilde{x}_+^i) - F_k(\tilde{x}_-^i)\right)\right)\right| \\
&\leq \frac{1}{4}\sum_{k=1}^K \sup_{x^i,\tilde{x}^i} \left|\left(F_k(x_+^i) - F_k(x_-^i)\right) - \left(F_k(\tilde{x}_+^i) - F_k(\tilde{x}_-^i)\right)\right| \\
&= \sum_{k=1}^K D_{ij}^{F_k}.
\end{align*}
Since the influence and derivative matrices always have non-negative entries the claim follows.
\end{proof}
We now bound the derivative matrix for some simple measures. The first calculation, for rank-1 quadratic forms, is well known. In the second result, we extend the calculation to quartic forms.

\begin{Lem} \label{lem:vectordob}
Let $F:\DC \to \mathbb{R}$, $F(x) = \langle u,x\rangle^2$, for some $u \in \mathbb{R}^n$. Then,

\begin{equation*}
\|D^F\|_{{\op}} \leq 2\|u\|^2_2.
\end{equation*}
\end{Lem}

\begin{proof}
A calculation shows that $\partial_i \partial_j F =2 u_i u_j$. Thus, we have $D^F_{ij} = 2|u_i u_j|$. This identity shows that if $\bar{u} = (|u_1|, \dots, |u_n|)$, then $D^F = 2\bar{u} \otimes \bar{u}$. Therefore, 

\begin{equation*}
\|D^F\|_{{\op}} =2 \|\bar{u}\|^2_2 = 2\|u\|^2_2.
\end{equation*}
\end{proof}

\begin{Lem} \label{lem:vectorproductdob}
Let $F:\DC \to \mathbb{R}$, $F(x) = \langle u,x\rangle^2 \langle v,x\rangle^2$, for some $u,v \in \mathbb{R}^n$. Then,

\begin{equation*}
\|D^F\|_{{\op}} \leq 12n \|u\|^2_2 \|v\|^2_2.
\end{equation*}
\end{Lem}
\begin{proof}
Fix $i,j \in [n]$ and write $F_1(x) = \langle u,x\rangle^2$ and $F_2(x) = \langle v,x\rangle^2$. Some straightforward calculations show
\begin{align*}
F_1(x) &\leq n \|u\|^2_2, \\
\partial_i F_1(x) &= 2u_i \langle u,x^i\rangle \leq 2\sqrt{n}|u_i|\|u\|_2, \\
\partial_i \partial_j F_1(x) &= 2u_i u_j.
\end{align*} 
Applying the Leibniz rule for discrete derivatives, we find
\begin{align*}
D^F_{ij} &\leq \max_x |\partial_i \partial_j F_1(x)| \max_x |F_2(x)| + \max_x |\partial_i \partial_j F_2(x)| \max_x |F_1(x)| \\
&\quad + \max_x |\partial_i F_1(x)| \max_x |\partial_j F_2(x)| + \max_x |\partial_i F_2(x)| \max_x |\partial_j F_1(x)| \\
&\leq 2n \left(|u_i u_j|\|v\|^2_2 + |v_i v_j|\|u\|^2_2 + 2|u_i v_j|\|u\|_2\|v\|_2 + 2|v_i u_j|\|u\|_2\|v\|_2\right).
\end{align*}
The bound follows by estimating the norm of the four rank-1 matrices, as in Lemma \ref{lem:vectordob}.
\end{proof}

We require one more result that explains how the spectral gap behaves with respect to general non-negative measures.
\begin{Lem}\label{lem:nonnegativesg}
	Let $\nu$ be a non-negative measure on $\DC$ and set $\tilde{\nu} = \frac{1}{\EE_\nu\left[1\right]}\nu$. Suppose that $C_{\poi}(\tilde{\nu})<\infty$. Then, for every $\vphi:\DC \to \RR$,
	$$\overline{\var}_{\nu}(\vphi) \leq C_{\poi}(\tilde{\nu})\mathcal{E}_{\nu}(\vphi),$$
	where $\mathcal{E}_{\nu}(\vphi) :=\EE_\nu\left[1\right]\mathcal{E}_{\tilde{\nu}}(\vphi)$.
\end{Lem}
\begin{proof}
	Since $\overline{\var}_{\nu}(\vphi) = \EE_{\nu}[\vphi^2] - \frac{\EE_{\nu}[\vphi]^2}{\EE_{\nu}[1]} = \EE_{\nu}[1]\left(\EE_{\tilde{\nu}}[\vphi^2]-\EE_{\tilde{\nu}}[\vphi]^2\right)= \EE_{\nu}[1]\var_{\tilde{\nu}}(\vphi)$, we get
	$$\overline{\var}_{\nu}(\vphi) \leq C_{\poi}(\tilde{\nu})\EE_{\nu}[1]\mathcal{E}_{\tilde{\nu}}(\vphi)= C_{\poi}(\tilde{\nu})\mathcal{E}_{\nu}(\vphi).$$
\end{proof}
 \subsection{Spectral gaps of tensor Ising models}
With the above estimates and our decomposition theorem, we can now provide a spectral gap estimate for general tensor Ising models.
	\begin{proof}[Proof of Theorem \ref{thm:tensorgap}]
		Fix a test function $\varphi:\DC\to \RR$. We apply Theorem \ref{thm:maindecompo} and obtain a decomposition
        \begin{equation*}
            \mu=\int \mu_{u,v,w,\ell,\psi} \eta (\dd u,\dd v,\dd w,\dd\ell, \dd\psi)\,,
        \end{equation*}
        with
\begin{equation*}
	 \mu_{u,v,w,\ell,\psi}(x) \propto \exp\left(\langle u, x \rangle^2 \langle v, x \rangle^2 + \langle x, w \rangle^2 + \langle \ell,x\rangle +\psi(x)\right),
	\end{equation*}
	and the following almost sure inequalities,
 	 $$\lVert u \lVert^2,  \lVert v \lVert ^2 \leq 2\sqrt{\lVert T\lVert_\inj},\ \  \lVert w\lVert^2\leq 4n\lVert T \lVert_\inj,\text{ and }\max\limits_{x\in \DC}|\phi(x)|\leq \delta,$$
 	 for some arbitrary $\delta>0$. 
	 Write 
    $$F_{u,v} := \langle u,x\rangle^2\langle v,x\rangle^2, \ G_w(x)= \langle w,x\rangle^2, \text{ and } L_\ell(x) = \langle \ell,x\rangle.$$ By Lemma \ref{lem:vectordob} and Lemma \ref{lem:vectorproductdob}, we get 
    \begin{align*}
       \lVert D^{F_{u,v}}\lVert_\op,&  \leq 48n\lVert T \lVert_\inj \\
        \lVert D^{G_{w_1}}\lVert_\op,& \leq 8n\lVert T\lVert_\inj\\
        \lVert D^{L_{\ell}}\lVert_\op &= 0,
    \end{align*} 
    where the last identity follows since $\partial_i\partial_j L_\ell(x) = 0$, for every $i,j \in [n]$. Moreover, note that for every $i,j\in[n]$ $|\partial_j\psi(x)| \leq 2\delta$ and $|\partial_i\partial_i\psi(x)|\leq 4\delta$, so $\|D^\psi\|_{\mathrm{op}} \leq n\delta.$ By
    applying Lemma \ref{lem:additivedob} to the each measure $\mu_{u, v, w,\ell,\psi}$, we get,
		\begin{align*}
		    \|A^{\mu_{u, v, w,\ell,\psi}}\|_{\op} &\leq  \lVert D^{F_{u,v}}\lVert_\op + \lVert D^{G_{w}}\lVert_\op + \lVert D^{L_{\ell}}\lVert_\op  + \lVert D^{\psi}\lVert_\op\\
            & \leq 48n \lVert T \lVert_\inj+ 8n \lVert T \lVert_\inj +\delta n=56n\lVert T \lVert_\inj + \delta n.            
            \end{align*}
         Hence, as long as $\delta$ is small enough, for the normalized measure $\widetilde{\mu}_{u, v, w,\ell,\psi}=\frac{\mu_{u, v, w,\ell,\psi}}{\EE_{\mu_{u, v, w,\ell,\psi}}\left[1\right]}$, by the Dobrushin condition \eqref{eq:dobru}, 
		$$C_\poi(\widetilde{\mu}_{\bar u,\bar v,\bar w,\ell})\leq\frac{1}{1-\sg n\|T\|_{\mathrm{inj}}-\delta n}.$$ Applying the variance condition of Theorem \ref{thm:maindecompo}, along with Lemma \ref{lem:nonnegativesg}, yields the following Poincar\'e inequality
		\begin{align*}
			\var_\mu(\varphi)&\leq\int \overline\var_{\mu_{u,v, w,\ell,\psi}}(\varphi)\eta(\dd u,\dd v, \dd w,\dd \ell,\dd \psi) +\delta
			\\
			&\leq \frac{1}{1-\sg n\|T\|_{\mathrm{inj}}-\delta n}\int\mathcal{E}_{\mu_{u,v, w,\ell,\psi}}(\varphi)\eta(\dd u,\dd v, \dd w,\dd \ell,\dd \psi) +\delta.
		\end{align*}
        To finish the proof, we shall use the fact that, even when considering general non-negative measures, the Dirichlet form is concave (see Lemma \ref{lem:dirichlet} in the appendix) to recompose into the Dirichlet form of $\mu$. Specifically, concavity affords the following inequality
        $$\int\mathcal{E}_{\mu_{u,v, w,\ell,\psi}}(\varphi)\eta(\dd u,\dd v, \dd w,\dd \ell,\dd \psi) \leq \mathcal{E}_\mu(\varphi).$$
        Together with the previous inequality, and since $\delta$ is arbitrary, we conclude
        $$	\var_\mu(\varphi) \leq \frac{1}{1-\sg n\|T\|_{\mathrm{inj}}} \mathcal{E}_\mu(\varphi).$$
	\end{proof}
    As a direct consequence, we now derive a spectral gap estimate for degree-$4$ spin glasses.
	\begin{proof}[Proof of Corollary \ref{cor:spinglass}]
	First, by thinking about $T$ as an operator on matrices, we consider its positive definite shifted version
	 $$\tilde T = T + \|T\|_{\mathrm{inj}}\mathrm{I}_{n^2}.$$
 Clearly $\tilde T$ is positive definite and $\|\tilde T\|_{\mathrm{inj}} \leq 2\|T\|_{\mathrm{inj}}.$
	 	Since $\mu$ is a measure on $\DC$, we also have the identity,
	 $$\mu(x) \propto\left(\beta \tilde T(x)\right).$$
	By Theorem \ref{thm:ben} we know that, for every $\varepsilon>0$, $\PP \left( \big\lvert n \lVert T\lVert_\inj-E_0(4)\big\lvert \leq \varepsilon \right)\xrightarrow{n\to\infty} 0 $. Therefore we get
    \begin{equation*}
        \lVert \Tilde{T} \lVert_\inj \leq 2\lVert T \lVert _\inj \leq (1 + o(1)) \frac{2E_0(4)}{n}.
    \end{equation*}
	It follows by Theorem \ref{thm:tensorgap} that for large enough $n$,
	$$C_{\poi}(\mu)\leq \frac{1}{1-\sg n\beta\| \Tilde{T}\|_{\inj}} \leq \frac{1}{1-2\cdot \sg\beta E_0(4)}.$$
    It remains to notice that by Theorem \ref{thm:ben}, $2\cdot \sg\cdot  E_0(4) \approx \spinsg$.
\end{proof}
\section{Spectral gap of tensor Curie-Weiss model} \label{sec:CW}

In this section, we focus on one specific tensor model, the Curie-Weiss model. Let us fix $p >1$ and consider the measure 
$$\mu_{p,\beta}(x) \propto \exp\left(\frac{\beta}{n^{p-1}}|\langle x, \mathbf{1}\rangle|^p\right)=\exp\left(\frac{\beta}{n^{p-1}}\left|\sum x_i\right|^p\right), \ \ \ \ x \in \{-1,1\}^n.$$
Let $\{X_t\}_{\{t \geq 0\}}$ be the iterates of Glauber dynamics. Since we wish to show that for certain temperatures $\beta$, the mixing time $X_t$ is slow, it will be enough to study the magnetization chain $S_t = \langle X_t, \mathbf{1}\rangle$. The idea is, that as observable of the process, the mixing time of $X_t$ is always bounded from below by the mixing time of $S_t$. Moreover, $S_t$ is a one-dimensional birth-and-death chain, and the form of the Curie-Weiss model affords a succinct analytic representation for its transitions. Using this representation, we shall thus prove.
\begin{Prop}\label{prop:TensorCW}
    Fix $p>1$ and set $\beta^* = \min\limits_{x > 0} \frac{\tanh^{-1}(x)}{px^{p-1}}$. Suppose that $\beta > \beta^*$ then there exists a constant $\alpha > 0$, 
 such that the mixing time of Glauber dynamics, associated to the Curie-Weiss measure $\mu_{p,\beta}$, is bounded from below by $e^{n^\alpha}$.
\end{Prop}  
Lemma \ref{lem:4CW} follows from Proposition \ref{prop:TensorCW} with a numerical calculation showing that for $p = 4$, $\beta_* \simeq 0.50425$

Since we are interested in lower bounds, it will be natural to work here in discrete time, and bound from below the number of required updates. Thus, to prove Proposition \ref{prop:TensorCW}, we shall first look at the transition matrix of $S_t$. For that, note that $S_t \in \ZZ \cap [-n,n]$ almost surely, and that $S_{t+1}-S_t \in \{-2,0,2\}$, for every $t \geq 0$. Let us understand each possible update separately. 

To understand the case $S_{t+1} = S_t + 2$, consider the set $N := \{i \in [n] | X_{t,i} = -1\}$, and observe that if $S_t = s$, necessarily $|N| = \frac{1-\frac{s}{n}}{2}n.$ Moreover for $S_{t+1}$ to be larger than $S_t$ it is necessary that GLD chose a coordinate from $N$. Let $I_t$ be the coordinate GLD chose at time $t$, so that
\begin{align*}
    \PP\left(S_{t+1} = S_t + 2|S_t = s\right) &= \PP\left(I_t \in N\right) \PP\left(S_{t+1} = S_t + 2|S_t = s, I_t \in N\right) \\
    &=  \frac{1-\frac{s}{n}}{2}\PP\left(S_{t+1} = S_t + 2|S_t = s, I_t \in N\right).
\end{align*}
Let us calculate the remaining probability term as
\begin{align*}
    \PP&\left(S_{t+1} = S_t + 2|S_t = s, I_t \in N\right) = \frac{e^{\frac{\beta}{n^{p-1}}(s+2)^p}}{e^{\frac{\beta}{n^{p-1}}(s+2)^p} + e^{\frac{\beta}{n^{p-1}}s^p}} \\
    &= \frac{e^{\frac{\beta}{n^{p-1}}\left((s+2)^p - s^p\right)}}{e^{\frac{\beta}{n^{p-1}}\left((s+2)^p - s^p\right)} + 1} = \frac{e^{\frac{\beta}{n^{p-1}}\frac{(s+2)^p - s^p}{2}}}{e^{\frac{\beta}{n^{p-1}}\frac{(s+2)^p - s^p}{2}} + e^{-\frac{\beta}{n^{p-1}}\frac{(s+2)^p - s^p}{2}}} \\
    &= \frac{1 + \tanh\left(\frac{\beta}{n^{p-1}}\frac{(s+2)^p - s^p}{2}\right)}{2}:=p^+_s.
\end{align*}
and we arrive at $\PP\left(S_{t+1} = S_t + 2|S_t = s\right) = \frac{1-\frac{s}{n}}{2}p^+_s$.
Similarly, for the case $S_{t+1} = S_t - 2$ we define $P := \{i \in [n] | X_{t,i} = +1\}$ and use the fact $|P| = \frac{1+\frac{s}{n}}{2}n$ along with similar calculations to arrive at 
$$\PP\left(S_{t+1} = S_t + 2|S_t = s\right) = \frac{1+\frac{s}{n}}{2}p^-_s,$$ with $p^-_s$ given by,
$$p^-_s:= \frac{1 - \tanh\left(\frac{\beta}{n^{p-1}}\frac{(s-2)^p - s^p}{2}\right)}{2}.$$
Finally, with the above estimates we arrive at 
$$\PP\left(S_{t+1} = S_t|S_t = s\right) = 1 - \frac{1-\frac{s}{n}}{2}p^+_s - \frac{1+\frac{s}{n}}{2}p^-_s.$$

Denote the normalized total magnetization as $\tilde{s}:= \frac{s}{n}$. So, since $\tanh$ is $1$-Lipschitz and bounded, we get the representation,
\begin{align} \label{eq:probs}
    \PP\left(S_{t+1} = S_t + 2|S_t = s\right) &= \frac{1-\tilde{s}}{4}\left(1+ \tanh(p\beta \tilde{s}^{p-1}) + O\left(\frac{1}{n}\right)\right) \nonumber \\
    \PP\left(S_{t+1} = S_t - 2|S_t = s\right) &= \frac{1+\tilde{s}}{4}\left(1 - \tanh(p\beta \tilde{s}^{p-1}) + O\left(\frac{1}{n}\right)\right).
\end{align}

The transitions make it clear that the model is symmetric to switching signs, we may thus suppose, with no loss of generality, that $X_0 = \mathbf{1}$, and so $S_0= n$. With this assumption, we will find a regime of $\beta$, for which exponential time must pass before $S_t < 0$. Again, the symmetry of the model implies symmetry for the invariant distribution and such a result will imply a lower bound on the mixing time.
To show this, we will capitalize on the following elementary estimate about biased random walks, see \cite[(2.8), Chapter XIV.2]{feller1968intro}.
\begin{Lem} \label{lem:biased}
    Let $M_t$ be a random walk on $\ZZ$ such that $|M_{t+1}-M_t| \in \{0,1\}$ with
    $$\PP\left(M_{t+1} = M_t+1\right)=(1-q),\ \ \ \ \PP\left(M_{t+1} = M_t-1\right)=q.$$
    Suppose that $q > \frac{1}{2}$, then,
    $$\PP\left(\min_{t \geq 0}M_t \leq M_0 - m\right)\leq \left(\frac{1-q}{q}\right)^m.$$
\end{Lem}
The idea is now clear, to show that from the starting configuration $S_0 = n$ it takes exponential time for $S_t$ to cross $0$, it will be enough to find an interval $[a,b] \subset [0,n]$ of linear size, such that $S_t$ has a uniform positive drift. In that case, Lemma \ref{lem:biased} will show that $S_t$ must take exponential time to cross this interval.

\begin{proof}[Proof of Proposition \ref{prop:TensorCW}]
In light of Lemma \ref{lem:biased} we would like to understand for which $\beta$ and $\tilde{s}$,
\begin{equation} \label{eq:bias}
    \frac{\PP\left(S_{t+1} = S_t + 2|S_t = s\right)}{\PP\left(S_{t+1} = S_t - 2|S_t = s\right)} > 1.
\end{equation}
With the expression in \eqref{eq:probs}, for large enough $n$, the above inequality is satisfied when 
$$\frac{(1 - \tilde{s})\left(1+ \tanh(p\beta \tilde{s}^{p-1})\right)}{(1 + \tilde{s})\left(1- \tanh(p\beta \tilde{s}^{p-1})\right)}> 1,$$
where $\tilde{s} = \frac{s}{n}.$
With some algebraic manipulations, this inequality is equivalent to,
$$\tanh(\beta p \tilde{s}^{p-1}) > \tilde{s} \implies \beta > \frac{\tanh^{-1}(\tilde{s})}{p\tilde{s}^{p-1}}.$$
Now, suppose that $\beta > \beta^*$. Then, by definition of $\beta^*$, there exists $\tilde{s} > 0$, which satisfies the inequality. Moreover, since $\frac{\tanh^{-1}(x)}{px^{p-1}}$ is continuous on $(0,1)$, there exists an interval $[a,b]\subset (0,1)$ such that 
$\tanh(\beta p \tilde{s}^{p-1}) > \tilde{s}$ for every $\tilde{s} \in [a,b]$. Equivalently, \eqref{eq:bias} is satisfied when $s \in [a\cdot n, b\cdot n]$.
In particular, with another continuity argument, we may assume that for $s \in [a\cdot n, b\cdot n]$,
\begin{equation*}
    \frac{\PP\left(S_{t+1} = S_t + 2|S_t = s\right)}{\PP\left(S_{t+1} = S_t- 2|S_t = s\right)} > 1 +\delta,
\end{equation*}
for some $\delta > 0$.

Thus, there exists $q > \frac{1}{2}$, and a biased random walk $M_t$ with
     $$\PP\left(M_{t+1} = M_t+1\right)=1-q,\ \ \ \ \PP\left(M_{t+1} = M_t-1\right)=q,$$
     such that the magnetization chain $S_t$ stochastically dominates (a lazy version of) $M_t$, whenever $S_t \in [a\cdot n, b\cdot n]$. We will use this fact to show that when $S_0 = n$, $S_t$ is confined to the positive half-line for an exponential period of time.

Thus, suppose that for some $t_0\geq 0$, $S_{t_0} = b\cdot n,$
and define the stopping time 
$$\tau := \inf\{t\geq t_0| S_t \notin [a\cdot n, b\cdot n]\}.$$
Applying Lemma \ref{lem:biased} and using the stochastic domination shows
$$\PP\left(S_\tau < a\cdot n| S_{t_0} = b\cdot n\right) \leq \left(\frac{1-q}{q}\right)^{(b-a)n}.$$
For the starting configuration $S_0 = n$, by repeating the above argument for every $t$ in which $S_t = b\cdot n$, and applying a union bound, we see that there exists a constant $\alpha$, such that 
$$\PP\left(\exists t \leq e^{n^{\alpha}} | S_t< 0\right)\leq \PP\left(\exists t \leq e^{n^{\alpha}} | S_t< a\cdot n\right)\leq e^{-n^\alpha}.$$
To finish the argument, by symmetry, for the Curie-Weiss measure $\mu_{p,\beta}(s< 0) = \frac{1}{2} +o(1),$ which shows that the mixing time of $S_t$ is bounded from below by $e^{n^\alpha}$. Since $S_t = \langle X_t,\mathbf{1}\rangle$, we obtain the same bound for the mixing time of Glauber dynamics, $X_t$.
\end{proof}
\section{Decomposing higher-degree tensors} \label{sec:higher} In this section we discuss high-level ideas on how to extend the decompositions from Section \ref{sec:TSL} to higher-degree tensors. As will be seen, our guarantees deteriorate rapidly when the degree increases. Thus, we opt to keep the discussion informal and track the main quantities of interest.

In Theorem \ref{thm:maindecompo} we regarded $T$ as an $n^2\times n^2$ matrix and reduced the measure into measures involving $n\times n$ matrices, effectively halving the degree of the involved tensors. It is then natural to extend our decomposition to every degree of the form $2^p$, for some $p \in \NN$, in which case we can iterate the decomposition until reaching tensors of simple forms. We will begin the discussion by explaining the case of degree-$8$ tensor  (see Figure \ref{fig:schema4}) before moving to general powers of $2$.

We do mention that with a proper embedding, our approach also applies to arbitrary tensors. If $T$ is a tensor of \emph{even} degree $m$ such that $2^{p-1}< m<2^p$ we can always consider (the symmetrized version of) $\hat{T} = T\otimes\left(\frac{1}{n}\mathrm{I}_n\right)^{\otimes \frac{2^p-m}{2}}$ which is a tensor of degree $2^p$ and for every $x \in \DC$ $T(x) = \hat{T}(x)$ and
$\|\hat{T}\|_{\mathrm{inj}} = \left(\frac{1}{n}\right)^{\frac{
2p-m}{2}}\|T\|_{\mathrm{inj}}$. If instead $m$ is \emph{odd} we shall require a more elaborate embedding. Write $z$ for the concatenated vector $(x^{\otimes \frac{m-1}{2}}, x^{\otimes \frac{m+1}{2}})$. The $m$-degree tensor $T$ can be represented as a bilinear form in $z$ as follows
$$T(x) = \left\langle z,\left[
\begin{array}{c|c}
    0 & J \\ 
    \hline
    J^T & 0
\end{array}
\right]z\right\rangle,$$
where $J \in \mathrm{Mat}(n^{\frac{m+1}{2}}, n^\frac{m-1}{2})$ and the above matrix can be regarded as an even degree tensor.
\subsection{Decomposing tensors of degree \texorpdfstring{$8$}{lg}}
Let us consider the measure $\mu(x) \propto \exp(T_8(x))$, for some degree-$8$ tensor $T_8$. We consider $T_8$ as an element of $\mathrm{Mat}(n^4,n^4)$ and apply a similar reasoning to Theorem \ref{thm:firstdecomp}. In that case, we can decompose the measure into a mixture of measures of the form
$$\exp\left(T_4^{(1)}(x) \otimes T_4^{(2)}(x) +T_4^{(3)}(x)\right),$$
where $T_4^{(1)}, T_4^{(2)}$ and $T_4^{(3)}$ are degree-$4$ tensors, which satisfy $\lVert T_4^{(1)}\lVert_\inj,\lVert T_4^{(2)}\lVert_\inj \leq 2\sqrt{\lVert T_8\lVert_\inj} $ and $\lVert T_4^{(3)}\lVert_\inj\leq 4n^2\lVert T_8\lVert_\inj$. We can now iterate the application of Theorem \ref{thm:firstdecomp}, and decompose each $T_4^{(i)}$ separately. This results in a finer decomposition into  the following type of measures,
\begin{equation*}
    \exp\left(\left(\langle x, M^{(1)}x\rangle^2+\langle x, R^{(1)}x\rangle\right)\left(\langle x, M^{(2)}x\rangle^2+\langle x, R^{(2)}x\rangle\right) +\langle x, M^{(3)}x\rangle^2+\langle x, R^{(3)}x\rangle\right),
\end{equation*}
where 
now $\{M^{(i)},R^{(i)}\}_{i=1}^3$ are matrices and
 $$\lVert M^{(i)}\lVert_\op\leq 2\sqrt{\lVert T_4^{(i)}\lVert_\inj},\ \ \lVert R^{(i)}\lVert_\op\leq 4n\lVert T_4^{(i)}\lVert_\inj, \ \ i=1,2,3.$$
 We now apply Theorem \ref{thm:seconddecomp} for the terms $\langle x, M^{(i)}x\rangle^2 $ as well as for $\langle x, R^{(i)}x\rangle $ (to formally apply the theorem, we can multiply by $\langle x, \frac{1}{n}I_n,x\rangle \equiv 1$), in total 9 times. We get a final decomposition into elements of the form 
\begin{equation*}
    \bar\mu(x)\propto
    \exp\left(P_1(x)P_2(x)+P_3(x)\right),
\end{equation*}
where  
\begin{equation*}
    P_j(x) = \langle u_1^{(j)},x\rangle^2 \langle u_2^{(j)},x\rangle^2+\langle u^{(j)}_3,x\rangle^2.
\end{equation*}
\begin{figure}[ht!]
    \centering
    \resizebox{15cm}{!}{\begin{tikzpicture}[
    arrow/.style={-Stealth, thick},
    line/.style={draw, thick},
    dottedline/.style={draw, thick, dotted}
]

\node (A) at (0, 0) {\( T_8 \)};
\node (B0) at (-2, -2) {$T_4^{(1)}$};
\node (B1) at (-1, -2) {$\otimes$};
\node (B2) at (0, -2) {$T_4^{(2)}$};
\node (B3) at (1, -2) {$+$};
\node (B4) at (2, -2) {$T_4^{(3)}$};
\node (C0) at (-4, -4) {$\big(T_1^{\otimes 2}\otimes T_1^{\otimes 2}+T_1^{\otimes 2} \big)$};
\node (C1) at (-1, -4) {$\otimes$};
\node (C2) at (0, -4) {$T_4$};
\node (C3) at (1, -4) {$+$};
\node (C4) at (2, -4) {$T_4$};
\node (D0) at (-7, -6) {$\big(T_1^{\otimes 2}\otimes T_1^{\otimes 2}+T_1^{\otimes 2} \big)$};
\node (D1) at (-4.5, -6) {$\otimes$};
\node (D2) at (-2, -6) {$\big(T_1^{\otimes 2}\otimes T_1^{\otimes 2}+T_1^{\otimes 2} \big)$};
\node (D3) at (1, -6) {$+$};
\node (D4) at (2, -6) {$T_4$};
\node (E0) at (-7, -8) {$\big(T_1^{\otimes 2}\otimes T_1^{\otimes 2}+T_1^{\otimes 2} \big)$};
\node (E1) at (-4.5, -8) {$\otimes$};
\node (E2) at (-2, -8) {$\big(T_1^{\otimes 2}\otimes T_1^{\otimes 2}+T_1^{\otimes 2} \big)$};
\node (E3) at (1, -8) {$+$};
\node (E4) at (4, -8) {$\big(T_1^{\otimes 2}\otimes T_1^{\otimes 2}+T_1^{\otimes 2} \big)$};
\draw[thick] (A) -- (B2);
\draw[thick] (B0) -- (C0);
\draw[dotted, thick] (B2) -- (C2);
\draw[dotted, thick] (B4) -- (C4);
\draw[dotted, thick] (C0) -- (D0);
\draw[thick] (C2) -- (D2);
\draw[dotted, thick] (C4) -- (D4);
\draw[dotted, thick] (D0) -- (E0);
\draw[dotted, thick] (D2) -- (E2);
\draw[thick] (D4) -- (E4);
\end{tikzpicture}}
    \caption{Decomposition of a order $8$ tensor $T_8$.}
    \label{fig:schema4}
\end{figure}
Figure \ref{fig:schema4} illustrates the procedure we just described. \\
The norms are bounded as follows:
\begin{equation*}
    \lVert u_1^{(j)}\lVert^2,\lVert u_2^{(j)}\lVert^2\leq 2\sqrt{\lVert T_4^{(i)}\lVert_\inj}\leq 
\begin{cases}
    2^{\frac{3}{2}}\lVert T_8\lVert^ \frac{1}{4}_\inj \quad &\text{ for } j=1,2\\
    8n\sqrt{\lVert T_8\lVert_\inj} \quad&\text{ for } j=3
\end{cases}
\end{equation*}
and 
\begin{equation*}
    \lVert u_3^{(i)}\lVert ^2\leq 2n\sqrt{\lVert T_4^{(i)}\lVert_\inj}\cdot \lVert M^{(i)}\lVert_\op\leq \begin{cases}
        8n\sqrt{\lVert T_8\lVert_\inj} & \quad\text{ for } j=1,2\\
        16n^3\lVert T_8\lVert_\inj &\quad \text{ for } j=3
    \end{cases}
\end{equation*}
At this point, it will make sense to try to understand any possible implications for spectral gaps. For that, we require the following extension of Lemma \ref{lem:vectorproductdob}, which we shall prove at the end of the section.
\begin{Lem}\label{lem:extt}
    Let $(u^{(k)})_{k \in [K]}$ be a family of vectors. Define $F:\DC \rightarrow \RR$, $F(x)= \prod_{k=1}^K \langle u^{(k)}, x \rangle^2. $ Then
    \begin{equation*}
        \lVert D^F \lVert _\op \leq 2K(2K-1) n^{K-1} \prod_{k=1}^K \lVert u^{(k)}\lVert^2.
    \end{equation*}
\end{Lem}By combining Lemma \ref{lem:extt} and the above estimates on the norms,
it follows that $\lVert D^{P_1P_2}\lVert_\op \leq 2^{13}n^3\lVert T_8\lVert_\inj$ and $\lVert D^{P_3}\lVert_\op\leq 800 n^3 \lVert T_8\lVert$. Thus, by Lemma \ref{lem:additivedob} we can bound the influence matrix associated to $\bar{\mu}$ as
\begin{equation*}
    \lVert A^{\bar{\mu}}\lVert_\op\leq 8992n^3\lVert T_8\lVert_\inj \,.
\end{equation*}
This implies that $\mu$ satisfies a Poincar\'e inequality whenever $\lVert T_8\lVert_\inj < \frac{1}{8992n^3}$.
\begin{proof}[Proof of Lemma \ref{lem:extt}]
    Define $F_k(x)= \langle u^{(k)}, x\rangle^2.$ First note that by Cauchy-Schwartz, $F_k(x) \leq n \lVert u^{(k)}\lVert^2$. For the derivatives we have
    \begin{align*}
        \partial_i F_k(x) &= 2 \lvert u_i^{(k)} \lvert \langle u^{(k)},x\rangle \leq 2\sqrt{n}\lvert u_i^{(k)} \lvert \lVert u ^{(k)}\lVert  \\
        \partial_i \partial_j F_k(x) &= 2\lvert u_i^{(k)}\lvert \lvert u_j^{(k)}\lvert.
    \end{align*}
Applying Leibniz rule, we find 
\begin{align*}
    D_{ij}^F &\leq \sum_{k=1}^K \max_x \lvert \partial_i \partial_j F_k (x)\lvert \prod_{l \neq k} \max_x \lvert F_l(x)\lvert + \sum_{1\leq k \neq l \leq K} \max_x \lvert \partial_i F_l(x)\lvert \max_x \lvert \partial_j F_k(x) \lvert \prod_{\ell \neq k,l} \max_x \lvert F_\ell(x)\lvert \\
     &\leq 2n^{K-1} \sum_{k=1}^K \lvert u_i^{(k)}\lvert \lvert u_j^{(k)}\lvert \prod_{l\neq k} \lVert u^{(l)} \lVert^2 +4n^{K-1} \sum_{1\leq k \neq l \leq K} \lvert u_i^{(l)} \lvert \lVert u^{(l)} \lVert \lvert u_j^{(k)} \lvert \lVert u^{(k)} \lVert \prod_{\ell\neq k,l} \lVert u^{(\ell)}\lVert^2 \,.
\end{align*}
We write $\bar{u}^{(k)}=(\lvert u_1^{(k)}\lvert, \dots, \lvert u_n^{(k)} \lvert) $ and it follows that 
\begin{align*}
    D^F &\preceq 2n^{K-1} \sum_{k=1}^K \bar{u}^{(k)} \otimes  \bar{u}^{(k)} \prod_{l \neq k} \lVert u^{(l)}\lVert^2+4n^{K-1}\sum_{k\neq l} \bar{u}^{(k)}\otimes \bar{u}^{(l)} \lVert u^{(l)} \lVert \lVert u^{(k)} \lVert \prod_{\ell\neq k,l} \lVert u^{(\ell)}\lVert^2\,,
\end{align*}
which implies 
\begin{align*}
    \lVert D^F \lVert_\op \leq 2Kn^{K-1}\prod_{k=1}^K \lVert u^{(k)}\lVert ^2+4K(K-1)n^{K-1}\prod_{k=1}^K \lVert u^{(k)}\lVert ^2.
\end{align*}
\end{proof}

\subsection{Decomposing tensors of degree \texorpdfstring{$2^p$}{lg}} 
For arbitrary $p\in\mathbb{N}$, if $T_{2^p}$ is a tensor of degree $2^p$ we can repeat and iterate the splitting procedure from the previous section. Figure \ref{fig:schema16} illustrates this procedure for the case of a tensor of degree $16$. In order to understand possible implications for bounds on the spectral gap, we will only analyze an idealized setting and restrict our attention to the leading term of the decomposition. This term corresponds to the a rank-one tensor of degree $2^p$ and has the form 
\begin{figure}[ht!]
    \centering
    \resizebox{16cm}{!}{\begin{tikzpicture}[
  level distance=4cm,
  level 1/.style={sibling distance=20cm},
  level 2/.style={sibling distance=7cm},
  level 3/.style={sibling distance=1.7cm},
  every node/.style={circle, draw},
  plain/.style={draw=none, fill=none} 
]

\node {$T_{16}$} 
  child foreach \i in {1,2,3} { 
    node {$T_8^{(\i)}$}
      child foreach \j in {1,2,3} { 
        node {$T_4^{(\i,\j)}$}
          child foreach \k in {1,2,3} { 
            node {$T_1^{(\i,\j,\k)}$}
          }
      }
  };

\node[plain] (A1) at (-10,-4) {\Huge $\otimes$}; 
\node[plain] (A2) at (10,-4) {\Huge $+$}; 
\node[plain] (B1) at (-23.5,-8) {\Huge $\otimes$}; 
\node[plain] (B2) at (-16.5,-8) {\Huge $+$};
\node[plain] (B3) at (-3.5,-8) {\Huge $\otimes$}; 
\node[plain] (B4) at (3.5,-8) {\Huge $+$};
\node[plain] (B5) at (16.5,-8) {\Huge $\otimes$}; 
\node[plain] (B6) at (23.5,-8) {\Huge $+$};
\end{tikzpicture}}
    \caption{Decomposition of a tensor of degree 16, $ T_{16} $. First, it is decomposed into $ T_8^{(1)} \otimes T_8^{(2)} + T_8^{(3)} $. Each $ T_8^{(i)} $ is then further decomposed as described in the previous section. The final decomposition consists of rank-1 tensors $ T_1^{(i,j,k)} $ with $ i, j, k = 1, 2, 3 $.}

    \label{fig:schema16}
\end{figure}
 \begin{equation*}
     L(x)=\prod_{j=1}^{2^{p-1}}\langle u^{(j)},x\rangle^2 
 \end{equation*}
 The bound on the norms of \( \lVert u^{(j)} \rVert^2 \) can be derived recursively. Let \( n_{2^p} := \lVert T_{2^p} \rVert_\inj \). As discussed in previous sections, it is necessary to recenter the tensors to ensure they are positive definite. This results in the inequality $
\lVert T_{2^{p-1}} \rVert_\inj \leq 2 \sqrt{\lVert T_{2^p} \rVert_\inj} $ for $p > 2$.
Consequently, we obtain the recursion $n_{2^{p-1}} \leq 2 \sqrt{n_{2^p}}$, and we need to understand the value $n_1$ which satisfies $\|u^{(j)}\|^2\leq n_1$. Since \( n_{2^p} := \lVert T_{2^p} \rVert_\inj \), a calculation produces the following solution to the recursion
\begin{equation*}
    n_1 \leq 2^{\sum_{k=0}^{p-2} 2^{-k}} \lVert T_{2^p} \rVert_\inj^{\frac{1}{2^{p-1}}}.
\end{equation*}
Thus, we conclude that $\lVert u^{(j)} \rVert^2 \leq 4 \lVert T_{2^p} \rVert_\inj^{\frac{1}{2^{p-1}}}.$ By Lemma \ref{lem:extt}, it follows that
 \begin{align*}
     \lVert D^L\lVert_\op&\leq 2^p(2^p-1) n^{2^{p-1}-1} \prod_{j=1}^{2^{p-1}} \lVert u^{(j)}\lVert ^2 \leq q^22^qn^{\frac{q}{2}-1}\lVert T_{2^p}\lVert_\inj\,,
 \end{align*}
 where $q=2^p.$
Thus, even without accounting for the lower order terms in the decomposition, the best spectral gap guarantees we can hope for are $\frac{1}{1-q^2 2^{q}n^{\frac{q}{2}-1}\|T_{2^p}\|_{\mathrm{inj}}}$. When $\|T_{2^p}\|_{\mathrm{inj}}$ is small enough but fixed, the former bound is of the order $O\left(q^22^{q}n^{\frac{q}{2}-1}\|T_{2^p}\|_{\inj}\right)$. Compare this bound to \cite[Theorem 2]{AJKTV} which instead scales like $O\left(q^2n^{\frac{q}{2}-1}\|T_{2^p}\|_{\mathrm{inj}}\right)$, and offers better guarantees for higher degrees.

\printbibliography
\begin{appendix}
    \section{Components of TSL}
\subsection{Smoothed projections} \label{sec:projections}
the construction of the smoothed projection from Lemma \ref{fct:Cprops} is based on the result from \cite[Lemma 2]{eldan2022spectral}, which we now recall. For $\delta > 0$, set $h(x) = e^{\frac{-x^2}{2\delta}}$, and if $H \subset \RR^n$ is a subspace and $v\in \RR^n$ is a vector, define
\begin{equation} \label{eq:2dCtdef}
C(H,v) = \mathrm{Proj}_H - (1-h(\|v_H\|))\bar{v}_H\otimes \bar{v}_H,
\end{equation}
where $v_H = \mathrm{Proj}_H(v)$ and $\bar v_H = \frac{v_{H}}{\|v_H\|}.$
Now, according to \cite[Lemma 2]{eldan2022spectral}, $C(H,\cdot)$ is Lipschitz continuous, satisfies 
\begin{equation*} \label{eq:smoothfirst}
0\preceq C(H,v) \preceq \mathrm{I}_H,\ \mathrm{Image}\left(C(H,v)\right)\subset H, \text{ and } \mathrm{Tr}(C(H,v)) \geq \mathrm{dim}(H) - 1.
\end{equation*}
Moreover
\begin{equation*} \label{eq:approxproj}
\|C(H,v)v\|\leq \delta.
\end{equation*}
The above properties essentially constitute Lemma \ref{fct:Cprops}.
\subsection{Bounded processes} \label{sec:boundedproc}
Our next goal is to prove the existence of the drifts in \eqref{eq:boundedprocdef1} and \eqref{eq:proddrift}, leading to the proof of Proposition \ref{prop:smallv}.
We begin with a simplified case, which applies as long as the matrix $C_t$ remains invertible, or equivalently, as long as the subspaces $\mathrm{Image}(T_t),\mathrm{Image}(M_t)$, defined by \eqref{eq:Tdef}, \eqref{eq:Mt} do not degenerate. The general case will obtained by iterating this construction over the decreasing sequence of image subspaces.
\begin{Lem} \label{lem:boundedproc}
	Let $C_t$ be a locally-Lipschitz diffusion coefficient supported on invertible positive definite matrices such that almost surely 
	$C_t \preceq \mathrm{I}_n.$
	Then, for every $\delta > 0$, there exists an adapted process $v^\delta_t$ in $\RR^n$ such that if $X_t^\delta$ solves the SDE 
	$$\dd X_t^\delta = C_t\dd B_t + C_t^2v^\delta_t\dd t,\ \ \ \|X^\delta_0\|^2 < \delta,$$
	then, almost surely, $\|X^\delta_t\|^2 \leq \delta$, for every $t \geq 0$. Moreover, we also have that 
	$C_tv^\delta_t$ does not explode in finite time.
\end{Lem}
\begin{proof}
	Fix $\delta > 0$, and for simplicity suppress the dependence on $\delta$ in the notation. Consider the processes,
	\begin{align} \label{eq:vtdef}
	v_t :=  -4\frac{nX_t}{\delta-\|X_t\|^2} \implies \dd X_t = C_t\dd B_t - 4\frac{nC_t^2X_t}{\delta-\|X_t\|^2}\dd t,
	\end{align}
	Let $Y_t:=\|X_t\|^2$, our goal is to show that $Y_t$ does not reach $\delta$ in finite time. This will also show that SDE defined by \eqref{eq:vtdef} has a unique strong solution; it has locally Lipschitz coefficients which imply the existence of a solution until a possible exposition time, which does not occur in finite time. Towards our goal, we apply It\^o's formula and obtain,
	\begin{align*}
	\dd Y_t &= 2X_t\cdot C_t\dd B_t - 4\frac{n\|C_tX_t\|^2}{\delta-\|X_t\|^2}\dd t + \|C_t\|^2_{\HS}\dd t 
	\\&= 2X_t\cdot C_t\dd B_t + \frac{\delta\|C_t\|^2_{\HS} - \|C_t\|^2_{\HS}Y_t - 4n\|C_tX_t\|^2}{\delta-Y_t}\dd t.
	\end{align*}
	We would like to transform $Y_t$ into a local martingale by some monotone function $f:\RR_+ \to \RR$.	
We make the choice $f(x) = \frac{1}{\delta - x}$ and obtain through It\^o's lemma,
	\begin{align*}
	\dd f(Y_t) &= 2f'(Y_t)C_tX_t\cdot \dd B_t + \left(\frac{\delta\|C_t\|_{\HS}^2 - \|C_t\|_{\HS}^2Y_t - 4n\|C_tX_t\|^2}{(\delta-Y_t)^3} + 2\frac{\|C_tX_t\|^2}{(\delta-Y_t)^3}\right)\dd t\\
	&=2f'(Y_t)C_tX_t\cdot \dd B_t + \frac{(\delta-Y_t)\|C_t\|_{\HS}^2  - (4n-2)\|C_tX_t\|^2}{(\delta-Y_t)^3}\dd t.
	\end{align*}
	To control the drift term, we restrict attention to times where $C_t$ is not too small. Thus, for $m \in \mathbb{N}$, let $\tau_m:=\inf\{t > 0| C_t \preceq \frac{1}{m}\mathrm{I}_n\}$ so that for $t < \tau_m$, $\|C_tX_t\|^2 \geq \frac{1}{m}Y_t$, and since $C_t \succeq \mathrm{I}_n$
	$$\frac{(\delta-Y_t)\|C_t\|_{\HS}^2  - (4n-2)\|C_tX_t\|^2}{(\delta-Y_t)^3} \leq \frac{(\delta-Y_t)n  - \frac{2n}{m}Y_t}{(\delta-Y_t)^3},$$
	as long as $Y_t \leq \delta$.
	 This turns $Z_t:=f(Y_{t\wedge \tau_m})$ into a positive super-martingale in the regime.
	
	Now, suppose towards a contradiction that there exists some finite $t' <\tau_m$, such that $Y_{t'} = \delta$. Since $Y_{t}$ has continuous paths, we may assume, with no loss of generality, that $Y_0 > \frac{m}{m+1}\delta$. Moreover, for the stopping time $\tau := \inf\{t >0: Y_t \in \{\frac{m+1}{m}\delta,\delta\}\},$ by restarting the processes if necessary, we may assume $\tau = t'$ and $Y_{\tau}=\delta$.
	
	On the other hand, the above argument shows that the stopped process $Z_{t\wedge \tau}$ is a positive super-martingale. By the martingale convergence theorem, $Z_{t\wedge \tau}$ converges to a finite limit, which necessarily must equal $f^{-1}(Y_\tau)$.
	Since $\lim\limits_{x\to \delta} f(x) = \infty$, this would make $Z_{t\wedge \tau}$ converge to infinity, which is a contradiction.
	Thus, almost surely, $Y_{t\wedge \tau_m} < \delta$, for every $t > 0$.  Take now $m \to \infty$ to conclude the first part of the proof.
	
	For the second part, set $R_t = \frac{\|X_t\|^2}{\delta - \|X_t\|^2}.$ Since $\frac{x}{(\delta-x)}= \frac{1}{\delta-x}-1$, we can see that 
	$dR_t = df(Y_t).$ So, if $\frac{\|C_tX_t\|^2}{\delta-\|X_t\|^2} > 1,$	$$\frac{(\delta-Y_t)\|C_t\|_{\HS}^2  - (4n-2)\|C_tX_t\|^2}{(\delta-Y_t)^3} < \frac{n  - n}{(\delta-Y_t)^2}=0.$$
	Again we can conclude that $R_t$ behaves like a super-martingale in the regime $\frac{\|C_tX_t\|^2}{1-\|X_t\|^2} \geq 1.$ 
	Since $\frac{\|C_tX_t\|^2}{1-\|X_t\|^2} \leq R_t$, with the same arguments as above we can now show $\lim\sup \frac{\|C_tX_t\|^2}{1-\|X_t\|^2} < \infty$. The claim follows since $\|C_tv_t\| = 4n\frac{\|C_tX_t\|}{1-\|X_t\|^2}.$
\end{proof}

With Lemma \ref{lem:boundedproc} we now prove that the processes in \eqref{eq:boundedprocdef1} and \eqref{eq:proddrift} are well-defined and prove Proposition \ref{prop:smallv}. 
\begin{proof}[Proof of Proposition \ref{prop:smallv}]
	Lemma \ref{lem:boundedproc} provides a construction for the drift $v^\delta_t$, as long as the matrix $C_t$ from \eqref{eq:Cdef} is invertible. However, due to the construction of the localization process, we expect $\mathrm{dim}\left(\mathrm{Image}(C_t)\right)$ to decrease eventually, and hence $C_t$ will stop being invertible. Combined with the definition of $v_t^\delta$ in \eqref{eq:vtdef} we see that, a-priori, the process defined in \eqref{eq:tensor_sl} will have a unique solution until some stopping time $T$, at which point either $C_t$ becomes degenerate, or experiences a jump due to a blowup of the processes $V_t$ and $W_t$.
	
	On the other hand the definition of $C_t$ in \eqref{eq:2dCtdef} shows that $\mathrm{Image}(C_t)$ is constant in $[0,T]$. By a limiting procedure, we can consider the measure $\mu_T$ and consider the process \eqref{eq:2dCtdef} on $\mathrm{Image}(C_T)$. Since $\mathrm{Image}(C_t)$ remains locally constant, we can continue to iterate this scheme. 
\end{proof}
    \section{Concavity of the Dirichlet Form} 
Here we shall prove that the Dirichlet form of GLD is concave with respect to mixtures of measures. This is a well-known fact in the classical setting that follows from Jensen's inequality. Below we show that the proof extends to the case of mixtures of general non-negative measures.
Recall that if $\nu$ is a non-negative measure and $\tilde{\nu} = \frac{\nu}{\EE_{\nu}[1]}$ we defined 
$\mathcal{E}_{\nu}(\vphi) = \EE_{\nu}[1]\mathcal{E}_{\tilde\nu}(\vphi)$.
\label{sec:dirichlet}
	\begin{Lem} \label{lem:dirichlet}
        Let $\mu$ be a probability measure on $\{-1,1\}^n$ which admits a decomposition as a mixture,
        $$\mu = \int\limits \mu_\theta \dd \eta(\theta),$$
        where for each $\theta$, $\mu_\theta$ is a non-negative measure on $\{-1,1\}^n$, and where $\eta$ is some mixing probability measure. Then, for every $\vphi:\{-1,1\}^n\to \RR$
        $$\int\limits\mathcal{E}_{\mu_\theta}(\vphi)\dd \eta(\theta) \leq  \mathcal{E}_\mu(\vphi) ,$$
        where $\mathcal{E}_\mu$ is the associated Dirichlet form for Glauber dynamics.
        \end{Lem}
	\begin{proof}
	    Write $\mathcal{E}_\mu(\vphi)$ explicitly (see \cite[Lemma 9]{eldan2022spectral} for example) as 
        $$\mathcal{E}_\mu(\vphi) = \sum\limits_{{x,y}\in\{-1,1\}^n}(\vphi(x)-\vphi(y))^2\frac{\mu(x)\mu(y)}{\mu(x)+\mu(y)}.$$
        In particular, according to our definition, we have for every $\theta$,
        $$\mathcal{E}_{\mu_\theta}(\vphi) = \sum\limits_{{x,y}\in\{-1,1\}^n}(\vphi(x)-\vphi(y))^2\frac{\mu_\theta(x)\mu_\theta(y)}{\mu_\theta(x)+\mu_\theta(y)}.$$
        It is straightforward to verify that the function $(x,y) \to \frac{xy}{x+y}$ is concave on $[0,\infty)\times [0,\infty)$. 
        Hence, for any $x,y\in \{-1,1\}^n$, since $(\mu(x),\mu(y)) =  \int\limits (\mu_\theta(x),\mu_\theta(y)) \dd \eta(\theta)$, by Jensen's inequality we obtain
        $$\int\frac{\mu_\theta(x)\mu_\theta(y)}{\mu_\theta(x)+\mu_\theta(y)}\dd\eta(\theta) \leq \frac{\mu(x)\mu(y)}{\mu(x)+\mu(y)}.$$
        We can now identify $\mu$ as a vector in $\RR_+^{2^n}$, and take a linear combination of the above inequality with positive coefficients, to obtain
        $$\int\limits\mathcal{E}_{\mu_\theta}(\vphi)\dd \eta(\theta) \leq \mathcal{E}_\mu(\vphi) .$$
	\end{proof}
\end{appendix}
\end{document}